\documentclass[twocolumn]{autart}

\newcommand{\integ}{\mathbb{Z}}

\newcommand{\map}[3]{#1: #2 \rightarrow #3}

\newcommand*{\QEDB}{\hfill\ensuremath{\square}}
\newcommand*{\QEDBL}{\hfill\ensuremath{\blacksquare}}
\newcommand\oprocendsymbol{\hbox{$\square$}}
\newcommand\oprocend{\relax\ifmmode\else\unskip\hfill%
\fi\oprocendsymbol}

\newcommand{\real}{\mathbb{R}} 

\newcommand{\mc}{\mathcal}

\newcommand{\sbs}[2]{{#1}_{\textup{#2}}}

\newcommand{\norm}[1]{ |{#1}|}

\newcommand{\tsp}{\mathsf{T}} 
 
\newcommand{\inv}{{\negat 1}} 
\newcommand{\negat}{\scalebox{0.75}[.9]{\( - \)}}

\newtheorem{theorem}{Theorem}[section]

\newtheorem{assumption}{Assumption}

\newtheorem{problem}{Problem}

\newcommand{\blue}[1]{{#1}}
\newcommand{\green}[1]{{#1}}

\usepackage{hyperref}
\usepackage{color}
\usepackage{graphicx}     
\graphicspath{ {./IMG/} }
\usepackage{amsmath}
\usepackage{amssymb}
\usepackage{mathrsfs}
\usepackage{subfigure}
\usepackage{url}
\usepackage{booktabs}
\usepackage{array}
\usepackage[table]{xcolor}
\usepackage{mathtools} 		
\usepackage{epstopdf}
\usepackage{cite}
\usepackage{mathtools} 
\usepackage{picins}
\usepackage{lineno}
\journal{Automatica}

 \usepackage{epstopdf}
\begin{document}

\begin{frontmatter}




\title{ 
Online Optimization of Switched LTI Systems Using Continuous-Time and Hybrid Accelerated Gradient Flows
}
\vspace{-0.6cm}
%
\author{Gianluca Bianchin}, 
\author{Jorge I. Poveda}, 
\author{Emiliano Dall'Anese}


\address{Department of Electrical, Computer, and Energy Engineering, University of Colorado, Boulder, CO, 80309 USA.}  
\vspace{-0.7cm}


\maketitle

\begin{abstract}
This paper studies the design of feedback controllers to steer a switching linear time-invariant dynamical system towards the solution trajectory of a time-varying convex optimization problem. We propose two types of controllers: (i) a continuous controller inspired by the online gradient descent method, and (ii) a hybrid controller that can be interpreted as an online version of Nesterov's accelerated gradient method with restarts of the state variables.  By design, the controllers continuously steer the system towards the time-varying optimizer without requiring knowledge of exogenous disturbances affecting the system. For cost functions that are smooth and satisfy the Polyak-\L ojasiewicz inequality, we demonstrate that the online gradient-flow controller ensures uniform global exponential stability when the time scales of the system and controller are sufficiently separated and the switching signal of the system varies slowly on average. For cost functions that are strongly convex, we show that the hybrid accelerated controller outperforms the continuous gradient descent method. When the cost function is not strongly convex, we show that the the hybrid accelerated method guarantees global practical asymptotic stability. 
\end{abstract}
\end{frontmatter}

\section{Introduction}

\vspace{-0.2cm}
In this paper, we investigate the use of online optimization algorithms 
for the control of switching dynamical systems. We consider linear 
time-invariant (LTI) plants with state $x\in\mathbb{R}^n$, output 
$y\in\mathbb{R}^p$, and dynamics
%
\begin{align}
\label{eq:plantModel1}
\dot x &= A_\sigma x + B_\sigma u + E_\sigma \omega_t 
:= P_{\sigma}(x,u,\omega_t), \nonumber\\
y &= C x + D \omega_t   ~~~~~~~~~~~~~:= h(x,\omega_t),
\end{align}
%
where $\omega_t:\mathbb{R}_{\geq0}\to\mathbb{R}^q$ is an unknown exogenous disturbance, $u \in \mathbb{R}^m$ is the control input,  $\sigma:\mathbb{R}_{\geq0}\to\mathcal{S}$ is a piece-wise continuous switching signal taking values in the finite set $\mathcal{S}:=\{1,2,\ldots,S\}$, with $S\in\mathbb{N}$, and 
$A_\sigma, B_\sigma, E_\sigma, C, D$ are matrices of appropriate 
dimensions. 
\blue{The  goal is to steer the inputs and outputs 
of~\eqref{eq:plantModel1} towards the time-varying solutions of the problem:}
\begin{align}\label{eq:optimizationObjective1}
&\min_{u,x,y} \;\; \phi_u(u) + \phi_y(y), \nonumber\\
&~~\text{s.t.}~~P_{\sigma}(x,u,\omega_t)=0,~~~y = h(x,\omega_t),
\end{align}
where $\map{\phi_u}{\real^m}{\real}$ and $\map{\phi_y}{\real^p}{\real}$ 
are 
functions that embed performance metrics associated with the steady-state
input and output of the system, respectively. This problem can be 
interpreted as an equilibrium-selection problem, where the 
objective is to select at every time the equilibrium points  of \eqref{eq:plantModel1} that minimize the cost in \eqref{eq:optimizationObjective1}. 
This class of optimization problems has emerged in several engineering 
applications \cite{brunner2012feedback,AJ-ML-PPJVDB:09,MC-ED-AB:20,hauswirth2020timescale,menta2018stability,lawrence2018optimal,lawrence2018linear,zheng2019implicit}, including power systems \cite{menta2018stability,MC-ED-AB:20} and transportation systems \cite{GB-FP:19-tits}, where the time-variability of $\omega_t$ 
precludes the off-line solution of ~\eqref{eq:optimizationObjective1} for the purpose of real-time control.  

\vspace{-0.2cm}
Feedback-based optimizing controllers for 
\eqref{eq:plantModel1}-\eqref{eq:optimizationObjective1} were studied 
in~\cite{menta2018stability} when $\omega_t$ is constant. The authors 
considered low-gain gradient-flow controllers of the form:

\vspace{-0.6cm}
\begin{equation}\label{gradient_flows}
    \dot{u}=-\eta \nabla\Phi(u,y),
\end{equation}

\vspace{-0.25cm}
where $\Phi$ is a modified cost function and $\eta>0$ is a small tunable gain. This approach was extended to smooth nonlinear plants and controllers in~\cite{hauswirth2020timescale} using singular perturbation tools \cite[Ch. 11]{HKK:96}. 
Joint stabilization and regulation problems related to \eqref{eq:optimizationObjective1} were the focus of 
\cite{lawrence2018linear,lawrence2018optimal} for a class of smooth systems. For LTI systems under time-varying disturbances, problem \eqref{eq:optimizationObjective1} was addressed in \cite{MC-ED-AB:20} via integral quadratic constraints, providing conditions that guarantee exponential stability and bounded 
tracking errors. Similar time-varying settings for feedback-linearizable plants were considered in~\cite{zheng2019implicit}. 

\vspace{-0.3cm}
\blue{
Despite the above line of work, optimization-based controllers for
systems with switching dynamics have not been studied yet. These 
systems are prevalent in engineering applications where plants are 
characterized by multiple operating modes; these  include transportation 
networks, where multiple modes originate due to the switching nature of  
traffic lights, and power grids, whose dynamics have multiple modes 
due to switching hardware. For these systems, it remains an open 
question whether  optimization-based controllers can still be applied, 
and under what conditions on the switching signal it is possible to 
guarantee their convergence.}
In this work, we provide an answer to these questions by presenting new 
stability results for optimization-based controllers applied to 
switched LTI systems. By using Lyapunov-based tools for set-valued hybrid
dynamical systems (HDS) \cite{RG-RGS-ART:12}  and the notion of 
input-to-state stability (ISS), we show that an average dwell-time 
constraint~\cite{JPH-ASM:99} is sufficient to guarantee closed-loop 
stability, provided that the time scales of the plant and the controller are suitably separated. 

\vspace{-0.3cm}
\blue{
One of the well-known disadvantages of gradient-flow methods in the
convex optimization literature is that their rate of convergence is bounded by 
the fundamental limit $\mathcal{O}(1/t)$ \cite{Wibisono1}. Naturally, 
this limitation extends to cases where gradient flows are utilized for controlling dynamical systems under time scale separation, as in~\eqref{gradient_flows}.
A natural question to ask is whether \emph{accelerated} methods, such
as those studied in~\cite{WS-SB-EC:14}, can also be used as feedback 
controllers for dynamical systems. We address this question by studying controllers inspired by a family of ODEs with dynamic momentum  of the 
form:}

\vspace{-0.6cm}
\blue{
\begin{equation}\label{ODENesterov}
\ddot{u}+\frac{p+\dot \tau}{\tau} \dot{u} + p^2 \tau^{p-2} k \nabla \Phi(u,y)=0,    
\end{equation}

\vspace{-0.3cm}
where $\tau$ denotes time, $p \geq 2$, and $k>0$ is a gain. Systems of this form have recently  received attention due to their ability to optimize smooth convex cost functions at a rate of $\mc O(1/\tau^p)$ \cite{WS-SB-EC:14,Wibisono1}. While different versions 
of~\eqref{ODENesterov} have been explored for classical optimization 
problems (see \cite{JIP-NL:19,annaswamytunners,HighResolution2018,UribeHeavyBall}), the authors in  \cite{WS-SB-EC:14,Wibisono1} showed that, when $p=2$, system \eqref{ODENesterov} models a continuous-time approximation of Nesterov's accelerated gradient method. However, while existing results have established convergence 
of~\eqref{ODENesterov} to solve optimization problems without plants in 
the loop, guaranteeing its convergence in the presence of plant dynamics
is not trivial.
Indeed, as recently shown in \cite[Sec. IV.B]{hauswirth2020timescale} 
via numerical experiments, the interconnection 
between~\eqref{ODENesterov}  and a dynamical plant can result in
instabilities even when $k \rightarrow 0$.
This observation finds a theoretical explanation through~\cite{poveda2019inducing}, where the authors show that 
for~\eqref{ODENesterov} no (strict) Lyapunov function exists due to absence of uniformity in the convergence 
(see~\cite[Prop. 1]{poveda2019inducing} and \cite[Thm. 1]{PovedaACC20}).} This prevents the direct application of standard singular perturbation tools \cite{TeelNesicTAC,HKK:96} to establish closed-loop stability. To overcome these challenges, in this paper we introduce a new feedback controller that combines
the continuous-time dynamics \eqref{ODENesterov} with discrete-time 
periodic resets. We show that these hybrid controllers guarantee robust approximate tracking as well as acceleration properties. \blue{To the 
best of our knowledge, the results of this paper are the first that 
incorporate switching plants and hybrid controllers in online 
optimization.}

\vspace{-0.4cm}
\section{ Preliminaries and Problem Statement}
\label{sec_notation}

\vspace{-0.3cm}
Given a compact set $\mc A \subset \real^n$ and $x \in \real^n$, we 
define $|x|_{\mc A}:=\min_{y \in \mc A} \|y-x\|_2$. 
When $\mathcal{A} =\{0\}$, $|x|_{\mc A}=|x|$ denotes the norm of $x$.
Given $v \in \real^n$ and $w \in \real^m$,
we let $(v,w) \in \real^{n+m}$ denote their concatenation.
For a matrix $M \in \real^{n \times n}$, we let
$\bar \lambda (M)$ and $\underline \lambda(M)$ denote 
its largest and smallest eigenvalues, respectively. 

\vspace{-0.3cm}
\subsection{Set-Valued Hybrid Dynamical Systems}
\vspace{-0.3cm}

We use the framework of HDS to analyze switching systems and hybrid algorithms using a common mathematical framework. A HDS with state $\varphi \in \real^n$, is given by
\begin{align}\label{eq:hybridSystemDef}
\varphi &\in C,~  \dot \varphi \in F(\varphi), 
& \varphi &\in D,~  \varphi^+ \in G(\varphi), 
\end{align}
where $F:\mathbb{R}^n\rightrightarrows\mathbb{R}^n$ 
and $G:\mathbb{R}^n\rightrightarrows\mathbb{R}^n$ are the 
flow and jump maps, respectively, whereas $C \subset \real^n$ and $D \subset 
\real^n$ are the flow and jump sets, respectively.
System \eqref{eq:hybridSystemDef} 
generalizes continuous-time systems ($D=\varnothing$) and discrete-time systems ($C=\varnothing$). Solutions
to \eqref{eq:hybridSystemDef} are parametrized by two time indices: a 
continuous index $t\in\mathbb{R}_{\geq0}$ that increases continuously 
whenever the system flows in $C$ as
$\dot{\varphi}(t,j) := \frac{d}{dt}\varphi(t,j) \in F(\varphi(t,j))$; and a discrete index  $j\in\mathbb{Z}_{\geq0}$ 
that increases by one whenever the system jumps in $D$ as $\varphi^+:=\varphi(t,j+1) \in G(\varphi(t,j))$. 
Solutions to \eqref{eq:hybridSystemDef} are defined on hybrid 
time-domains~\cite[Def. 2.3]{RG-RGS-ART:12}, namely, subsets of  
$\mathbb{R}_{\geq0} \times\mathbb{Z}_{\geq 0}$ defined 
as the union of intervals  $[t_j,t_{j+1}]\times \{j\}$, with
$0 = t_0\leq t_1 \leq \ldots$, and where the last interval can be closed 
or open on the right. In compact form, we denote by $\text{dom}(\varphi)$
the domain of $\varphi$. 
%

\vspace{-0.2cm}
We study switching signals $\sigma$ that satisfy an 
\emph{average dwell-time} (ADT) condition \cite{JPH-ASM:99} of the form 
$N(t,s) \leq N_0 + \frac{t-s}{\tau_d}$, $\forall~0 \leq s \leq t$, where
$N(t,s)$ denotes the number of discontinuities of $\sigma$
in the interval $(s,t]$, and $\tau_d>0$ is called the dwell-time. The ADT 
condition guarantees that system \eqref{eq:plantModel1} has at most $N_0$
switches at any time, and finitely-many switches in any finite time interval. As shown in \cite[Ch. 2]{RG-RGS-ART:12}, HDS of the form \eqref{eq:hybridSystemDef} offer a  mathematical model to capture any signal $\sigma$ satisfying the ADT condition. In particular, every switching signal $\sigma$ satisfying ADT can be generated by a HDS with state $\chi = (\tau,\sigma)\in\mathbb{R}_{\geq0}\times\mathcal{S}$, 
data $(C_{\chi},F_{\chi},D_{\chi},G_{\chi})$ given by:
%
\begin{align}\label{eq:dwellTimeDynamics_hybrid}
C_{\chi}&:=[0, N_0]\times\mathcal{S},~~~~F_{\chi}(\chi):=[0,\tau_d^{-1}]\times \{0\}, \nonumber \\
D_{\chi}&:=[1, N_0]\times\mathcal{S},~~~~G_{\chi}(\chi):=\{\tau -1\}\times \mc S . 
\end{align}
%
Moreover, every signal $\sigma$ generated by the HDS \eqref{eq:dwellTimeDynamics_hybrid} has a hybrid time domain that satisfies the ADT condition. 

\vspace{-0.2cm}
\subsection{Problem Statement}

\vspace{-0.3cm}
We next formalize the problem of interest. For any
fixed $u \in \real^m$, $\sigma\in\mathcal{S}$, and $\omega_t \in \real^q$, we write the steady-state output of \eqref{eq:plantModel1} as:
%
$y_{\sigma}= G_{\sigma} u + 
H_{\sigma} \omega_t$,
%
where $G_{\sigma}:= -C A_{\sigma}^\inv B_{\sigma}$ and $H_{\sigma} := D -C A_{\sigma}^\inv E_{\sigma}$. We will impose the following assumptions.

\vspace{-0.2cm}
\begin{assumption}
\label{ass:stabilityPlant}
For each 
$\sigma \in \mc S$, and each symmetric matrix $P_{\sigma} \succ 0$, there exists a unique symmetric matrix $Q_{\sigma} \succ 0$, such that $A_{\sigma}^\tsp P_{\sigma} + P_{\sigma} A_{\sigma} = - Q_{\sigma}$.
\end{assumption}

\vspace{-0.2cm}
\begin{assumption}
\label{ass:sharedEquilibria}
\blue{For each $\omega_t \in\real^q$ and each
$u \in \real^m$, there exists a unique $x \in \real^n$, such that $A_{\sigma} x + B_{\sigma} u + E_{\sigma} \omega_t=0$, 
for all $\sigma  \in \mc S$.}
\end{assumption}

\vspace{-.2cm}

%
%
Under Assumption \ref{ass:stabilityPlant}, $A_\sigma$ is Hurwitz and
therefore invertible. On the other hand, Assumption \ref{ass:sharedEquilibria} is 
common for the analysis of switched 
systems~\cite{JPH-ASM:99,RG-RGS-ART:12} and it guarantees that all 
the modes have a common equilibrium (see also 
Remark~\ref{rem:average}), and that the input-output maps are common
across the modes, i.e., $G:=G_{\sigma}$ and $H:=H_{\sigma}$ for all 
$\sigma \in \mc S$. Under these assumptions, we can 
rewrite~\eqref{eq:optimizationObjective1} as:
\begin{align}\label{eq:optimizationObjectiveUnconstrained}
\min_{u} \;\;\; & \phi_t(u) := \phi_u(u) + \phi_y(G u + H \omega_t).
\end{align}
%
%
Note that every solution to \eqref{eq:optimizationObjective1} 
is a solution to~\eqref{eq:optimizationObjectiveUnconstrained},
however, the inverse implication holds only when 
$(A_{\sigma}, C_{\sigma})$, $\sigma \in \mc S$, is observable. 
Since we will focus 
on~\eqref{eq:optimizationObjectiveUnconstrained}, observability is 
not necessary in the subsequent analysis.
\blue{For simplicity, we  assume that for each $\omega_t$ 
problem~\eqref{eq:optimizationObjectiveUnconstrained} has a unique 
solution $u_t^*$, and that the mapping 
$\omega_t\mapsto u_t^*(\omega_t)$ is smooth, with $\omega_t$ 
satisfying the following assumption\footnote{\blue{In the next sections, and with some abuse of notation, we will retain the subscript $t$ to emphasize the dependence of $x^*_t$ and $u^*_t$ on the exogenous disturbance $\omega_t$. We will later drop this subscript in our stability analysis.}}.}

\vspace{-0.1cm}
\begin{assumption}\label{assumption_regularity}
\blue{The function $t\mapsto\omega_t$ is generated by an (unknown) Lipschitz continuous exosystem
\vspace{-.1cm}
\begin{equation}\label{exosystemdynamics}
\vspace{-.1cm}
\dot{\omega}_t=\Pi(\omega_t),~~~\omega_t\in \Gamma\subset\mathbb{R}^q,
\end{equation}
with $\Gamma$ being forward invariant and compact.}
\end{assumption}

\vspace{-0.2cm}
\blue{Assumption \ref{assumption_regularity} is standard for 
regulation problems with exogenous inputs or 
references~\cite{LM-AT:10,VG:07}, and it guarantees that 
$t\mapsto\omega_t$ is continuously differentiable and bounded.} 
%
%
%
\begin{rem} \label{rem:average}
{\normalfont 
When $G_{\sigma}$ and $H_{\sigma}$ are not common across modes, but 
all modes admit a common equilibrium point, we can 
define the average map $y = \sbs{G}{av} u + \sbs{H}{av} \omega_t$, where 
$\sbs{G}{av} := \sum_{\sigma \in \mc S} \alpha_{\sigma} G_{\sigma}$ and 
$\sbs{H}{av} := \sum_{\sigma \in \mc S} \alpha_{\sigma} H_{\sigma}$, 
with $0 \leq \alpha_{\sigma} \leq 1$, $\sum_{\sigma} \alpha_{\sigma}=1$. In this case, 
\eqref{eq:optimizationObjectiveUnconstrained} can be generalized to  
$\min_{u} \phi_u(u) + \phi_y(\sbs{G}{av} u + \sbs{H}{av} \omega_t)$. This scenario often emerges in transportation systems \cite{GB-FP:19-tits}.\QEDB
}
\end{rem}
%
\vspace{-.2cm}
In the remainder, we use $z:=(x,u)$ for the joint state of the plant
and the control signal, and $z_t^*:= (x_t^*, u_t^*)$ to denote the 
unique vector that satisfies the following equations for all times: 
$0 = A_\sigma x_t^* + B_\sigma u_t^* + E_\sigma \omega_t$, 
$\sigma\in \mathcal{S}$, and 
$0 = \nabla \phi_u(u_t^*) + \nabla \phi_y(G u_t^* + H\omega_t)$.  
In words, the components of $z_t^*$ correspond to the equilibria  
of~\eqref{eq:plantModel1} and the time-varying critical point 
of~\eqref{eq:optimizationObjectiveUnconstrained}. 
The problem focus  of this work is formalized next.
\vspace{-0.2cm}
\begin{problem}\label{prob:tracking}
Let $\tilde{z}:=z-z_t^*$ denote the tracking error. Design an 
output-feedback controller for \eqref{eq:plantModel1} 
such that for any unknown exogenous signal $t\mapsto \omega_t$, the tracking error converges to a neighborhood of the origin, whose size is parameterized by the time-variation of $\omega_t$, i.e., by $|\dot{\omega}_{t}|$.
\QEDB
\end{problem}

\vspace{-0.2cm}

\section{Main Results}
\label{sec:3}
\vspace{-.3cm}

To address Problem~\ref{prob:tracking}, in this paper we propose two different output-feedback
controllers: the first based on a gradient-descent flow, and the 
second based on a hybrid accelerated gradient method.
\begin{figure}[t!]
\centering
\includegraphics[width=.5\columnwidth]{%
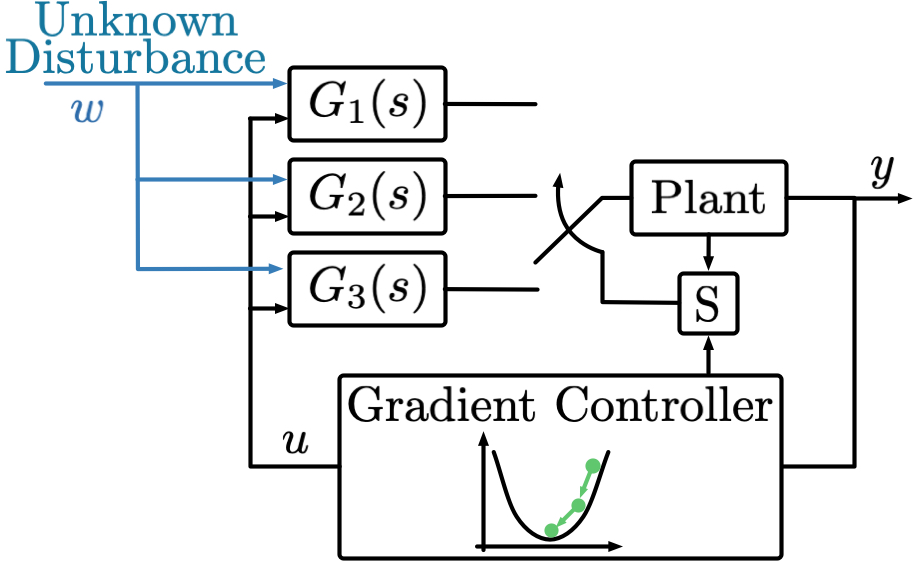}%
\includegraphics[width=.5\columnwidth]{%
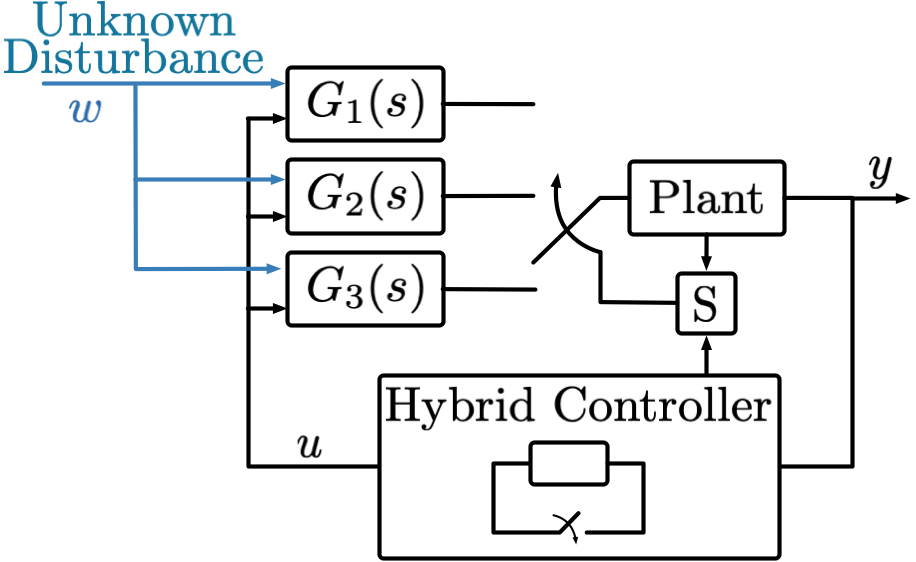}
\caption{\small{(left) A switched system in feedback with a gradient flow controller. (right) A switched system in feedback with a hybrid controller. ``S'' denotes a supervisory controller that actuates the switching between the modes of the plant.}}
\label{fig:plantControllerLoop} 
\end{figure} 

\vspace{-.2cm}
\subsection{Feedback Control via Online Gradient Descent}
\label{subsection_gradientflow}

\vspace{-0.3cm}
We first study the solution of the tracking 
Problem~\ref{prob:tracking} via gradient descent flows of the 
form~\eqref{gradient_flows}. The left scheme of
Fig.~\ref{fig:plantControllerLoop} illustrates the approach. For 
such systems, the following two assumptions are standard.

\vspace{-0.2cm}
\begin{assumption}\label{ass:lipschitzGradient}
The functions $\phi_u(\cdot)$ and $\phi_y(\cdot)$ are continuously 
differentiable, and their gradients are globally Lipschitz with 
constants $\ell_u>0$ and $\ell_y>0$, respectively.
\end{assumption}
\vspace{-0.2cm}
\begin{assumption}\label{ass:PLinequality}
For all $t\geq0$, the function $\phi_t(u)$: 
(a) is radially unbounded, and 
(b) satisfies the Polyak-\L ojasiewicz (PL) inequality, namely, 
$\norm{\nabla \phi_t(u)}^2 \geq 2 \mu ( \phi_t(u) - \phi_t(u_t^*))$,
for some $\mu > 0$ and for all $u \in \real^m$.
\end{assumption}
\vspace{-0.2cm}
Under Assumption  \ref{ass:lipschitzGradient}, the mapping
$u\mapsto \phi_t(u)$ has a globally Lipschitz gradient with 
Lipschitz constant $\ell := \ell_u + \ell_y \norm{G}^2$. Similarly, Assumption \ref{ass:PLinequality} implies that 
$\phi_t(u) - \phi_t(u^*) \geq \frac{\mu}{2} \norm{u - 
u_t^*}^2$, $\forall~u\in\mathbb{R}^m$. 
%
%
%
%

\vspace{-0.2cm}
To design the controller we note that if $\omega_t$ and $H$ were 
{known}, the following dynamics can be shown to converge to 
$u_t^*$ under Assumptions \ref{ass:lipschitzGradient} and 
\ref{ass:PLinequality} (see, for example,~\cite{absil2006stable}):
\vspace{-0.1cm}
\begin{align}
\label{eq:Fgs}
\hspace{-.2cm}
\dot u =F_{GS}(u,\omega_t):= - \nabla \phi_u(u) - G^\tsp \nabla \phi_y(G u +  H \omega_t).
\hspace{-.2cm}
\end{align}

\vspace{-0.6cm}
\blue{
When $\omega_t$ and $H$ are {unknown}, we propose to approximate the 
steady-state output $G u + H \omega_t$ with instantaneous 
feedback from the plant, leading to:}
\begin{align}\label{eq:gradient_controller}
\blue{\dot u = F_G(u,y,\sigma):=- \eta_\sigma ( \nabla \phi_u(u) + G^\tsp \nabla \phi_y(y)),}
\end{align}

\vspace{-0.4cm}
where $\eta_\sigma>0$ is a mode-dependent tunable gain. 

\vspace{-0.3cm}
The system obtained by interconnecting plant~\eqref{eq:plantModel1}, 
the switching signal generator~\eqref{eq:dwellTimeDynamics_hybrid}, 
and the controller~\eqref{eq:gradient_controller} leads to a HDS of 
the form \eqref{eq:hybridSystemDef}, denoted by $\mathcal{H}_G$, with
state $\varphi = (x, u,\chi,\omega_t)$, continuous-time~dynamics:
\blue{
\begin{align}
\label{eq:interconnectedGradFlow}
\dot x &= P_{\sigma}(x,u,\omega_t), &
\dot{u} &= {F_G(u,h(x,\omega_t), \sigma)}, \nonumber\\
\dot{\chi} &\in F_{\chi}(\chi),    &
\dot{\omega}_t &=\Pi(\omega_t)
\end{align}
}
flow set $C:=\real^n \times  \real^m \times C_{\chi}\times\Gamma$, discrete-time dynamics:
\begin{equation}
x^+=x,~~u^+=u,~~\chi^+\in G_{\chi}(\chi),~~\omega^+=\omega,   
\end{equation}
and jump set 
$D:=\real^n \times  \real^m \times D_{\chi}\times\Gamma$. 

Next, we provide a result that establishes an explicit tracking 
bound for \eqref{eq:interconnectedGradFlow}.
To this end, we 
require that the controller gain satisfies  $\eta_{\sigma} < \bar \eta_{\sigma}$, where 
\begin{align}
\label{eq:etaBoundGradient}
\bar \eta_{\sigma} =
\frac{(1-\kappa)^2}{2-\kappa}
\frac{ \underline{\lambda}(Q_{\sigma})}{
\ell_y \norm{C} \norm{G} \norm{P_{\sigma} A_{\sigma}^\inv B_{\sigma}} },
\end{align}
where $(P_{\sigma},Q_{\sigma})$ are as in 
Assumption~\ref{ass:stabilityPlant}, and $\kappa\in(0,1)$. Moreover, we define the following constants:

\vspace{-0.7cm}
\begin{small}
\begin{subequations}
\label{eq:aBaraUnderbarGradient}
\begin{align} 
\theta_{\sigma} &:= \green{\frac{1}{
1 + 2 \norm{P_{\sigma} A_{\sigma}^\inv B_{\sigma}}}},\\
\bar a_\sigma &:=\frac{1}{\eta} \max \left\{(1\negat\theta_{\sigma})\frac{\ell}{2},\theta_{\sigma} \bar
\lambda(P_{\sigma}) \right\},\label{barasigma}\\
\underline a_{\sigma}  &:=\frac{1}{\eta} \min \left\{ (1\negat\theta_{\sigma})\frac{\mu}{2},  
\theta_{\sigma} \underline \lambda(P_{\sigma}) \right\},\label{underbarasigma}\\
b_{\sigma} &=  \frac{\kappa}{2} \min \left\{ 
2\mu \eta_\sigma , 
\frac{\underline \lambda (Q_{\sigma})}{\bar \lambda (P_{\sigma})}
\right\},\label{bsigma}\\
d_{\sigma} &=\frac{2}{\kappa} \max \left\{
\frac{\ell_y \norm{H}\norm{G}}{\eta_\sigma \mu^2},
\frac{2 \norm{P_{\sigma} A_{\sigma}^\inv E_{\sigma}} }{\underline \lambda(Q_{\sigma})}
\right\},
\end{align}
\end{subequations}
\end{small}\noindent
where $(\mu,\ell,\ell_y,\ell_u)$ are as in 
Assumptions~\ref{ass:lipschitzGradient} and \ref{ass:PLinequality}.
\begin{theorem}\label{thm:EISS-Gradient}
Suppose that Assumptions \ref{ass:stabilityPlant}-\ref{ass:PLinequality} hold. If 
$\eta_{\sigma} \in (0, \bar \eta_{\sigma})$ for all 
$\sigma \in \mc S$ and the dwell-time satisfies 
$ \tau_d > \frac{\ln a}{\min_\sigma b_\sigma}$, 
then \blue{the tracking error $\tilde z = z -z_t^*$
of the system $\mathcal{H}_G$ satisfies:
\begin{align}\label{ISS_bound}
\hspace{-.3cm} |\tilde z(t,j)|
\leq  a_0 e^{-\frac{b_0 t + c_0 j}{2}} 
\vert \tilde z(t,j) \vert 
+ a_0d_0 \sup_{0 \leq \tau \leq t} 
\norm{\dot \omega_\tau}, \hspace{-.3cm}  
\end{align}
}\noindent

\vspace{-0.6cm}
where $a_0 = e^{\frac{\varrho}{2} N_0} \sqrt{a}$, 
$b_0 = \min_{\sigma} b_{\sigma} - \frac{\varrho}{\tau_d}$,
$c_0 = \varrho - \ln a$, 
$d_0 = \max_{\sigma} d_{\sigma}$, 
with
$a := \frac{\max_{\sigma} \bar a_{\sigma}}{
\min_{\sigma}\underline a_{\sigma}}$ and $\varrho>0$ is any constant 
that satisfies 
$\ln a < \varrho < \tau_d \min_{\sigma} b_{\sigma}$. 
\QEDB
\end{theorem}
%

\blue{
This result establishes that a sufficiently-small gain 
$\eta_{\sigma}$ guarantees exponential convergence of $z$ to a 
neighborhood of the optimal trajectory $z_{t}^*$. 
As characterized by~\eqref{eq:etaBoundGradient}, the upper
bound on the controller gain is  proportional to the rate of 
convergence of the open-loop plant 
$\underline \lambda(Q_\sigma)/\bar \lambda(P_\sigma)$,
and inversely proportional to the Lipschitz constant of the cost 
function $\ell_y$. Moreover, as characterized by \eqref{bsigma}, the 
rate of decay of the tracking error is governed by the minimum 
between the rate of convergence of the controller (namely, 
$2 \mu \eta_\sigma$) and the rate of convergence of the plant 
(namely, $\underline \lambda(Q_\sigma)/\bar \lambda(P_\sigma)$).
}

\vspace{-0.2cm}
\begin{rem}
\normalfont 
The constants $(a_0,b_0,c_0,d_0)$ characterized in Theorem \ref{thm:EISS-Gradient}
depend on: (i)~
$(\mu,\ell)$, which characterize the smoothness and gradient dominance of the \emph{steady state} cost function \eqref{eq:optimizationObjectiveUnconstrained}; (ii) the 
matrices $(A_{\sigma},B_{\sigma},C_{\sigma},D_{\sigma},E_{\sigma})$ and $(Q_{\sigma},P_{\sigma})$, which govern the dynamics of the switched plant \eqref{eq:plantModel1}; and (iii) $(\tau_d,N_0)$ that govern the
 behavior of the signal $\sigma$ in \eqref{eq:dwellTimeDynamics_hybrid}.
\QEDB
\end{rem}

\vspace{-0.2cm}
\begin{rem}
{\normalfont 
Assumption \ref{ass:PLinequality} is fundamental to guarantee 
exponential convergence. 
Indeed, as $\mu \rightarrow 0$, we have that $b_{\sigma} \rightarrow 0 $ and 
$d_{\sigma} \rightarrow \infty$. Similar scenarios were recently investigated in \cite{menta2018stability,hauswirth2020timescale}. In contrast to these results, Theorem \ref{thm:EISS-Gradient} accounts for \emph{time-varying} disturbances, \emph{switching} plant dynamics, and establishes an explicit \emph{exponential} bound, as opposed to asymptotic convergence.
}
\QEDB
\end{rem}

\vspace{-0.1cm}
\subsection{Feedback Control via Hybrid Gradient Descent}
\label{sec:4}
\vspace{-0.2cm}

We now address Problem \ref{prob:tracking} by proposing a feedback 
controller inspired by the accelerated gradient 
method~\eqref{ODENesterov}. 
\blue{
To design the controller, we adapt \eqref{ODENesterov} as follows: 
first we rewrite~\eqref{ODENesterov} as a set of first-order ODEs by 
letting $p=2$ and by defining the variables $u_1:=u$ 
and $u_2:=\frac{\tau}{2}\dot{u}_1+u_1$; second, we introduce 
an auxiliary state $u_3$ that models the evolution of a timer and 
replaces the temporal variable $\tau$, thus  leading to the 
dynamics:}
\vspace{-0.2cm}
\begin{subequations}\label{naiveflow}
\begin{align}
\dot{u}_1&=\frac{2}{u_3}(u_2-u_1)\\
\dot u_2&= - 2 k u_3\left( \nabla \phi_u(u_1) 
+ G^\tsp \nabla \phi_y(G u_1 +  H \omega_t)\right),\\
\dot{u}_3&=\frac{1}{2}.
\end{align}
\end{subequations}

\vspace{-.4cm}
\blue{

The choice of the variables 
$(u_1,u_2,u_3)$ is inspired by accelerated momentum-based optimization and 
estimation algorithms proposed in e.g. \cite[Eq. (14)]{Wibisono1}, 
\cite{annaswamytunners}, and \cite{JIP-NL:19}. Further, we note 
that the choice $\dot{u}_3=1/2$, also used in \cite{JIP-NL:19} to 
solve \emph{static} optimization problems, is motivated by our 
Lyapunov-based analysis, presented below.

\vspace{-.2cm}
\begin{rem}\label{rem:lackUniformity}
As shown in \cite[Prop. 1]{poveda2019inducing} and 
\cite[Thm. 1]{PovedaACC20}, the convergence of \eqref{naiveflow} 
lacks uniformity with respect to the initial value of $u_3$. In turn,
this  precludes the application of standard multi-time scale techniques using quadratic-type Lyapunov functions \cite[pp. 453]{HKK:96} or regular perturbations techniques, see \cite[Thm. 1]{TeelNesicTAC}.
\QEDB
\end{rem}
\vspace{-.2cm}
}

\blue{
Similarly to \eqref{eq:Fgs}, the dynamics \eqref{naiveflow} require 
knowledge of $H$ and $\omega_t$ to be implemented. Hence, we propose 
to approximate the steady-state map $Gu_1+H\omega_t$ with 
the instantaneous  output $y$ of the dynamical system. 
While this modification leads to an accelerated version 
of the controller proposed in Section~\ref{subsection_gradientflow}, 
empirical and theoretical evidence suggests that such modifications 
are not sufficient to guarantee tracking of the optimal trajectories 
(see Remark~\ref{rem:lackUniformity}).
To address this limitation, we introduce discrete-time {resets} of the state
variables of \eqref{naiveflow}, which resemble the ``restarting'' 
heuristics used in the literature of machine 
learning~\cite{WS-SB-EC:14,Wibisono1,o2015adaptive} and 
hybrid control \cite{Prieur}.} The proposed controller is then hybrid, with
state $u=(u_1,u_2,u_3)$, continuous-time dynamics

\vspace{-0.7cm}
\begin{small}
\begin{equation*}
\dot{u}=F_H(u,y,\sigma):=\eta_\sigma \left(\begin{array}{c}
 \dfrac{2}{u_3} (u_2-u_1)\\
 - 2k u_3  \left( \nabla \phi_u(u_1) +  G^\tsp \nabla \phi_y(y) \right) \\ 
\frac{1}{2} 
\end{array}\right),
\end{equation*}
\end{small}

\vspace{-0.5cm}
where \blue{$\eta_{\sigma}>0$ is a mode-dependent tunable gain}; 
flow set $C_H=\mathbb{R}^n\times\mathbb{R}^n\times[\delta,\Delta]$, 
where $\Delta>\delta>0$ are tunable parameters characterizing the 
restarts of the timer variable; discrete-time dynamics
\begin{equation*}
u^+=G_H(u):=R_0u,~~R_0=\left[\begin{array}{ccc}
I  & 0 & 0\\
r_0I & (1-r_0)I & 0\\
0 & 0 & \frac{\delta}{\Delta}
\end{array}\right],
\end{equation*}
where $r_0\in\{0,1\}$ is the \emph{reset policy} of the controller; 
and jump set $D_H=\mathbb{R}^n\times\mathbb{R}^n\times\{\Delta\}$. 

\vspace{-.3cm}
\blue{
In words, the logic of the controller is as follows: when the timer 
$u_3$ is equal to $\Delta$, the controller states $(u_1, u_2)$ are 
re-initialized and the timer variable is reset to the value $\delta$.
Notice that when $r_0=0$, only the timer $u_3$ is reset to $\delta$, 
whereas when $r_0=1$ both the momentum variable $u_2$ is 
re-initialized to $u_1$ and the timer $u_3$ is reset to $\delta$. 
By recalling the definition of $u_2$, we notice that the latter 
restarting policy corresponds to setting the momentum $\dot{u}_1$ to zero.
}

\vspace{-.3cm}
The interconnection of the hybrid controller and the switched LTI plant \eqref{eq:plantModel1} leads to a HDS with state $\varphi=(x,u,\chi,\omega_t)$, \blue{continuous-time dynamics:
\begin{align}\label{flowshybrid}
\dot x &= P_{\sigma}(x,u,\omega_t), &
\dot{u} &= {F_H(u,h(x,\omega_t), \sigma)}, \nonumber\\
\dot{\chi} &\in F_{\chi}(\chi),    &
\dot{\omega}_t &=\Pi(\omega_t)
\end{align}
and flow set $C=\mathbb{R}^n\times C_H\times C_{\chi} \times \Gamma$.}
Jumps in the closed-loop system can be triggered by both switches of 
the plant and by resets of the controller.  Therefore, the 
discrete-time dynamics are captured by the inclusion:
\begin{equation}\label{jumpshybrid}
\varphi\in D:=D_1\cup D_2,~~~\varphi^+\in G(\varphi),
\end{equation}
where the set-valued map $G$ and the sets $D_1,D_2$ are
\begin{equation}\label{definition_set_jumps}
G(\varphi):=\left\{\begin{array}{l}
G_1(\varphi),~~~~\varphi\in D_1:=\mathbb{R}^n\times C_H\times D_{\chi}\times\Gamma\\
G_2(\varphi),~~~~\varphi\in D_2:=\mathbb{R}^n\times D_H\times C_{\chi}\times\Gamma\\
G_1(\varphi)\cup G_2(\varphi),~~~\varphi\in D_1\cap D_2,
\end{array}\right.
\end{equation}
\blue{with $G_1(\varphi):=\{x\}\times\{u\}\times G_{\chi}(\chi) \times \{\omega_t\}$ and 
$G_2(\varphi):=\{x\}\times\{G_H(u)\}\times \{\chi\} \times \{\omega_t\}$.}
We denote the hybrid 
system~\eqref{flowshybrid}-\eqref{definition_set_jumps} in compact form 
by $\mathcal{H}_H$.

\blue{
\vspace{-.2cm}
\begin{rem}
Notice that the model 
\eqref{flowshybrid}-\eqref{definition_set_jumps} naturally captures 
non-uniqueness of solutions that can emerge when the plant and the 
controller are in their jump sets simultaneously. 
In these scenarios, arbitrarily-small disturbances may force the 
plant or the controller to jump before the other, and a well-posed 
model of the hybrid dynamics allows us to capture both behaviors of 
the system as the disturbance vanishes. 
\QEDB
\end{rem}
}

\vspace{-0.1cm}
Next, we provide error-tracking guarantees for 
$\mc H_H$. To this end, we first consider the case where the disturbance 
$\omega_t$ is constant, and we show that the (time-invariant) optimizer 
is globally practically {asymptotically} stable.

\vspace{-0.1cm}
\blue{
\begin{rem}
Since we first focus on asymptotic stability properties, it suffices to 
consider plants \eqref{eq:plantModel1} with a single mode
(i.e., $\mathcal{S}:=\{\sigma\} $). Indeed, if each individual mode 
in $\mathcal{S}$ leads to asymptotic stability of the 
closed-loop system, then semi-global practical asymptotic stability 
(with respect to $\tau_d^{-1}$ 
in~\eqref{eq:dwellTimeDynamics_hybrid}) for the system with multiple 
modes will follow directly from \cite[Corollary 7.28]{RG-RGS-ART:12}.
\QEDB
\end{rem}
}
\vspace{-0.2cm}

We begin by introducing an inverse Lipschitz-type assumption.

\vspace{-0.2cm}
\begin{assumption} \label{ass:reverseLipschitz}
The function $u\mapsto \phi_t(u)$ is convex, \blue{radially 
unbounded}, and for each $\nu_0 >0$ there exists $\ell_0>0$ such that
$\norm{u-u_t^*} \geq  \nu_0$ for some $u \in \real^m$ implies
$\norm{\nabla \phi_t(u)} \geq \ell_0 \norm{u- u_t^*}$.   
\QEDB
\end{assumption}

\vspace{-0.2cm}
Next, we require that the controller gain satisfies 
$\eta_\sigma < \bar \eta_\sigma$, where

\vspace{-0.9cm}
\begin{small}
\begin{equation*}
    \bar{\eta}_H=
\frac{\underline \lambda(Q_{\sigma})}{2 \ell_y \norm{C}\norm{G}}
\min \left\{
\frac{1}{2 k \delta \Delta \norm{P_{\sigma} A_{\sigma}^\inv B_{\sigma}}},
\frac{\theta\ell_0^2 \delta \kappa^2 k}{2 (1-\theta) \ell \Delta \norm{C}\norm{G}}
\right\}.
\end{equation*}
\end{small}
with $\kappa\in(0,1)$ and

\vspace{-0.8cm}
\begin{small}
\begin{equation*}
\theta = \frac{\ell_y k \Delta \norm{C}\norm{G}}{\ell_y k \Delta  \norm{C}\norm{G} + 2 \delta \norm{P_{\sigma} A_{\sigma}^\inv B_{\sigma}}}.
\end{equation*}
\end{small}

\vspace{-0.4cm}
Using this gain we obtain the following result:

\vspace{-0.2cm}
\begin{theorem}\label{thm:EISS-Nesterov-Practical}
\blue{
Let Assumptions~\ref{ass:stabilityPlant}-\ref{ass:lipschitzGradient} 
and \ref{ass:reverseLipschitz} hold, and assume that
$\omega_t:= \omega\in\mathbb{R}^q$ is constant, 
that the plant has a single mode $\mc S =\{ \sigma\}$, and let
the reset policy be $r_0=0$.
}
If $\eta_{\sigma} \in (0,\bar \eta_\sigma)$, then any solution of 
$\mc H_H$ satisfies:

\vspace{-0.7cm}
\begin{small}
\begin{align*}
&a)~ \limsup_{t+j \rightarrow \infty}|z(t,j)-z^*|
\leq \nu_0,~~~~~~~~~\forall(t,j)\in\text{dom}(\varphi).\\
&b)~ \phi(u_1(t,j))-\phi(u^*)\leq \frac{\alpha_j}{u_3^2(t,j)}+\nu_0,~\forall(t,j)\in\text{dom}(\varphi),t>\underline{t}_j,
\end{align*}
\end{small}\noindent 

\vspace{-0.4cm}
where $\alpha_j>0$, $z:=(x,u_1)$, $z^*:=(x^*,u^*)$,
and $\underline{t}_j:=\min\{t\geq0:(t,j)\in\text{dom}(\varphi)\}$.
\QEDB
\end{theorem}

The convergence result of Theorem \ref{thm:EISS-Nesterov-Practical} is  global, but of ``practical'' nature. Namely, convergence is achieved only to a $\nu_0$-neighborhood of the optimal set via the choice of $\bar{\eta}_\sigma$, which, in turn, depends on the constant $\ell_0$. 
Two main comments are in order. First, the upper bound on the 
controller gain $\bar{\eta}_\sigma$ shrinks to zero when either 
$\delta\to0^+$ or $\Delta\to\infty$. 
Since $\delta=0$ and $\Delta=\infty$ correspond to the (non-restarted 
ODE) \eqref{ODENesterov}, our analysis suggests that for 
\eqref{ODENesterov} there may exist no gain $\eta>0$ that guarantees 
stability for the closed-loop system: a similar observation was recently 
made in \cite[Sec. IV.B]{hauswirth2020timescale}. Second, the result establishes that the error in the \emph{steady-state} cost function decreases (outside a neighborhood of the optimal point) at a rate of order $\mathcal{O}(c_j/u_3^2)$ during the $j^{th}$ interval of flow, where $c_j$ is constant in each interval. Thus, the larger the difference $\Delta-\delta$, the larger the size of the intervals where this bound holds. This can be seen as a semi-acceleration property that holds during flows. Indeed, using the definition of $u_3$, during the first interval of flow we have that $u_3(t,0)=u_3(0,0)+0.5t$, and the error in the cost decreases at a rate $\mathcal{O}(c_0/t^2)$; see also \cite{WS-SB-EC:14,Wibisono1,JIP-NL:19} for similar bounds in static optimization problems. We also note that, as revealed by the proof (presented below), and in contrast to the case of 
optimization problems without plants in the loop, the quantity $\alpha_j$ 
in item \textit{b)} explicitly depends on the LTI plant 
\eqref{eq:plantModel1} via the matrix $P_{\sigma}$. 
Finally, we note that as $\nu_0 \rightarrow 0$, the controller gain might satisfy $\eta \rightarrow 0$, as shown in the 
following example.

\vspace{-0.1cm}
\begin{exmp}
\label{eq:reverserLipschitz}
Let $\phi(u_1) = \frac{1}{\theta}(\norm{u_1-u^*}^\theta + \norm{Gu_1+H\omega_t}^\theta)$, with $\theta \in \integ_{\geq 2}$, which is not strongly convex when $\theta>2$. Also, let $G=[g_{ij}]$ with 
$g_{ij}\geq0$ for all $i,j$, as it is the case in compartmental systems \cite{GB-FP:19-tits}. It then follows that
$\norm{\nabla \phi(u_1)} \geq \norm{u_1-u_t^*}^{\theta-1} \geq \nu_0^{\theta-2} 
\norm{u_1-u_t^*}$. Thus, $\phi_t$ satisfies Assumption  
\ref{ass:reverseLipschitz} with 
$\ell_0 = \nu_0^{\theta-2}$. Note that in this case, as $\nu_0\to0^+$, the admissible constant $\ell_0$ and the controller gain $\bar{\eta}_\sigma$ shrink to zero.
This relation suggests that as the size of the neighborhood satisfies
$\nu_0 \rightarrow 0$, the controller gain also satisfies 
$\eta \rightarrow 0$.
\QEDB
\end{exmp}

\vspace{-0.1cm}
Next, under the following strong convexity assumption, we will show that 
the hybrid controller can solve Problem 1 with an exponential rate of
convergence. 

\vspace{-0.2cm}
\begin{assumption}
\label{ass:strongConvexity}
There exists $\mu \in \real_{>0}$ such that 
$\phi_t(u) \geq \phi_t(u') + \nabla \phi_t(u')^\tsp (u-u') + \frac{\mu}{2} 
\norm{u-u'}^2$ holds for all $ u, u' \in \real^m$, and all $t\geq0$.
\end{assumption}
\vspace{-0.2cm}

\blue{
To establish an exponential tracking bound, we require that the 
controller gain satisfies $\eta_\sigma < \bar \eta_\sigma$, where

\vspace{-0.7cm}
\begin{small}
\begin{align}
\label{eq:etaBoundNesterov}
\bar{\eta}_\sigma &= \green{ 
\frac{(1-\kappa)^2}{16}
\frac{\delta \underline \lambda(Q_\sigma)}{
k \Delta \ell_y \norm{C}\norm{G} \norm{P_\sigma A_\sigma^\inv B_\sigma}},} &
\kappa & \in (0,1).
\end{align}
\end{small}\noindent 
Moreover, we define the following constants:

\vspace{-0.6cm}
\begin{small}
\begin{subequations}
\begin{align}
\label{eq:coefficientsNesterov-a}
\theta_{\sigma} &= \green{\frac{k \Delta\ell_y \norm{C}\norm{G}}{
 k \Delta\ell_y \norm{C}\norm{G} + 2  \delta^\inv \norm{P_\sigma A_\sigma^\inv B_\sigma}}},\\
\label{eq:coefficientsNesterov-b}
\bar a_{\sigma} &:= 
\green{\max \left\{ (1-\theta_{\sigma}) \frac{1 + k \ell \Delta^2}{2},
\theta_{\sigma} \bar \lambda(P_{\sigma}) \right\}},\\
\label{eq:coefficientsNesterov-c}
\underline a_{\sigma} &:= 
\green{\min \left\{(1-\theta_{\sigma}) \frac{1 + 2 k \ell \Delta^2}{4}, 
\theta_{\sigma}	\underline \lambda(P_{\sigma}) \right\}},\\
\label{eq:coefficientsNesterov-d}
b_{\sigma} &= \green{\frac{\kappa}{2} 
\min \left\{ \frac{2\eta_\sigma}{\Delta(1+2 k \ell \Delta^2)}, 
\frac{\eta_\sigma  k \delta \mu}{2(1+2 k \ell \Delta^2)},
\frac{\underline \lambda(Q_{\sigma})}{\bar \lambda (P_{\sigma})}
\right\}}, \hspace{-.2cm} \\
\label{eq:coefficientsNesterov-e}
d_{\sigma} &= \green{ 
\frac{2}{\kappa} \max \left\{ 
\frac{
\bar{d}_{\sigma}
}{
\eta_\sigma \min \{\Delta^\inv, k \delta \mu/4 \}},
\frac{2 \norm{P_{\sigma} A_{\sigma}^\inv E_{\sigma}}}{\underline \lambda(Q_{\sigma})}
\right\}}.
\end{align}
\end{subequations}
\end{small}

\vspace{-0.5cm}
\noindent
where $\bar{d}_{\sigma}:= \sqrt{2} k \Delta^2 \ell_y \norm{H} \norm{G} 
+ \sqrt{2} \max_{\omega \in \Gamma} \norm{\frac{\partial u^*}{\partial \omega}}(k \Delta^2 + \frac{1}{2}(1 + 2 k \Delta^2 \ell))$. 
}

\begin{theorem}
\label{thm:EISS-Nesterov}
Let Assumptions \ref{ass:stabilityPlant}-\ref{ass:lipschitzGradient}
and \ref{ass:strongConvexity} hold and let the reset policy be 
$r_0=1$. \blue{
If $\eta_{\sigma} \in (0, \bar \eta_{\sigma})$ for all 
$\sigma \in \mc S$, the timer thresholds satisfy
$\Delta^2 - \delta^2 >\frac{1}{2k \mu}$, and the dwell-time satisfies
$\tau_d > \frac{\ln a}{\min_\sigma b_{\sigma}}$, then the tracking error $\tilde z = z - z_t^*$ of the system $\mc H_H$  satisfies:
\begin{align}\label{ISS_bound}
\hspace{-.3cm} |\tilde z(t,j)|
\leq  a_0 e^{-\frac{b_0 t + c_0 j}{2}} 
\vert \tilde z(t,j) \vert 
+ a_0d_0 \sup_{0 \leq \tau \leq t} 
\norm{\dot \omega_\tau}, \hspace{-.3cm}  
\end{align}
where $a_0 = e^{\frac{\varrho}{2} N_0} \sqrt{a}$, 
$b_0 = \min_{\sigma} b_{\sigma} - \frac{\varrho}{\tau_d}$,
$c_0 = \varrho - \ln a$, 
$d_0 = \max_{\sigma} d_{\sigma}$, 
with
$a := \frac{\max_{\sigma} \bar a_{\sigma}}{
\min_{\sigma}\underline a_{\sigma}}$ and $\varrho>0$ is any constant 
that satisfies 
$\ln a < \varrho < \tau_d \min_{\sigma} b_{\sigma}$. 
}
\QEDB
\end{theorem}

\vspace{-0.25cm}
The result of Theorem  \ref{thm:EISS-Nesterov} requires three 
types of conditions on the controller parameters:
(i) a sufficiently-small choice of the controller gain $\eta_\sigma$,
(ii) a quadratic-like dwell-time condition 
$\Delta^2 - \delta^2 >1/2k\mu$, which gives a lower bound on the 
restarting frequency, and (iii) an average dwell-time condition 
imposed on the switching signal $\sigma$. 
\begin{rem}\label{final_remark}
The result of Theorem \ref{thm:EISS-Nesterov} leverages the resets of the momentum $\dot{u}$. Similar ``restarting'' techniques have been studied in literature of optimization and machine learning \cite{WS-SB-EC:14,Wibisono1}, with optimal restarting frequencies presented in \cite[Thm. 3.1]{o2015adaptive}, \cite[Sec. 3.2.1]{JIP-NL:19}. For example, in \cite[Sec. 3.2.1]{JIP-NL:19}, it is shown that, as $\delta\to0^+$,  the choice $\Delta=e\sqrt{\frac{1}{2k\mu}+\delta^2}$ with $k=1/2\ell$ can achieve exponential convergence of order $\mathcal{O}(e^{-t\sqrt{\mu/\ell}})$. This is particularly advantageous in problems with condition numbers satisfying $\ell/\mu\gg1$. In the context of feedback-based optimization, the advantages of using momentum are recovered as the time scale separation between the plant and the controller increases. Numerical experiments are presented in Section \ref{sec:6}.
\QEDB
\end{rem}
\blue{We conclude by noting that the closed-loop systems $\mc H_H$ with constant disturbances $\omega_t=\omega$ are robust to \emph{small bounded additive} disturbances $t\mapsto e(t)$ acting on the states and dynamics of the system. Indeed, by construction, the hybrid system $\mc H_H$ is well-posed and it satisfies the Basic Conditions \cite[Assump. 6.5]{RG-RGS-ART:12}. 
Therefore, the robustness property follows directly by \cite[Cor. 7.27]{RG-RGS-ART:12}.}
\vspace{-0.3cm}
\section{Proofs}
\label{sec_proofs}
\vspace{-0.3cm}
\blue{
To prove our main results, we use tools from HDS theory \cite{RG-RGS-ART:12}. We model the closed-loop system as a time-invariant and well-posed HDS described by the interconnection between the exosystem \eqref{exosystemdynamics}, the plant \eqref{eq:plantModel1}, the switching generator \eqref{eq:dwellTimeDynamics_hybrid}, and the controller. For this interconnection, we construct a hybrid Lyapunov function to guarantee stability of the closed-loop dynamics.
The construction of the Lyapunov functions is inspired by singular perturbation arguments \cite{HKK:96}, adjusted to account for the switching dynamics and the hybrid controllers. This construction allows us to derive conditions on $\eta_\sigma$ and on the dwell-time parameters of $\sigma$ to guarantee suitable asymptotic or exponential input-to-state stability with respect to $\dot{\omega}$. Since the closed-loop system can be modeled as a time-invariant HDS, in our analysis we drop the subscript $t$ from the signal $\omega$, the optimizer $u^*(\omega)$, and the cost function \eqref{eq:optimizationObjectiveUnconstrained}, which we write as an explicit function of $u$ and $\omega$ of the form $\phi(u,\omega)$.
} 
\vspace{-0.3cm}
\subsection{Proof of Theorem \ref{thm:EISS-Gradient}}
\vspace{-0.3cm}

We divide the proof into four main steps.

\vspace{-0.1cm}
\blue{
\textbf{Step 1.} Consider the change of variable $\tilde{u}:=u-u^*(\omega)$, and $\tilde x := x -\left(-A_{\sigma}^\inv B_{\sigma} (\tilde{u}+u^*(\omega)) - A_{\sigma}^\inv E_{\sigma} \omega\right)$, which denote the tracking error of the plant and the controller, respectively. The dynamics of $\tilde{u}$ are:
}

\vspace{-0.8cm}
\begin{small}
\begin{align}
\dot{\tilde{u}}&=\dot{u}-\dot{\overbrace{u^*(\omega)}}=
- \eta\Big( \nabla \phi_u(\tilde{u}+u^*(\omega)) \nonumber\\
&\quad + G^\tsp \nabla \phi_y\big(C\tilde{x}+G (\tilde{u}+u^*(\omega))
+H\omega\big)\Big)-\dot{\overbrace{u^*(\omega)}}\notag\\
&=:\eta\psi_c(\tilde x,\tilde{u}, \omega)-\dot{\overbrace{u^*(\omega)}},\label{errorplantdynamics}
\end{align}
\end{small}
\vspace{-0.8cm}

and note that $\psi_c(0,\tilde{u},\omega)=-\nabla_u \phi(\tilde{u}+u^*(\omega),\omega)$. Also

\vspace{-0.8cm}
\begin{small}
\begin{align}
\dot{\tilde{x}}&=
A_{\sigma}\tilde{x}+\eta A_{\sigma}^{-1}B_{\sigma}\psi_c(\tilde x,\tilde{u}, \omega)+A_{\sigma}^{-1}E_{\sigma}\dot{\omega}.\label{errorcontroldynamics}
\end{align}
\end{small}
\vspace{-0.3cm}
Next, we consider the Lyapunov-like function

\vspace{-0.5cm}
\begin{small}
\begin{align}
\label{eq:LyapunovFunction_Gradient}
 V_{\sigma}(\tilde{z},\omega) = \frac{(1-\green{\theta_\sigma})}{\eta} V_1(\tilde{u},\omega) + \frac{\green{\theta_\sigma}}{\eta} V_{2}(\tilde {x},\sigma),
\end{align}
\end{small}
\vspace{-.8cm}

with $\tilde{z}=(\tilde{x},\tilde{u})$ and $\theta_{\sigma}\in(0,1)$ \green{is as in \eqref{eq:aBaraUnderbarGradient}}, and
$V_1(\tilde{u},\omega) := \phi(\tilde{u}+u^*(\omega),\omega) - \phi(u^*(\omega),\omega)$,
$V_2(\tilde x,\sigma) := \tilde x^\tsp P_{\sigma} \tilde x, 
$
with $P_{\sigma}$ given by Assumption \ref{ass:stabilityPlant}. The function $\green{V_\sigma}$ is a convex combination of a Lyapunov-like function for the plant and a Lyapunov-like function for a standard gradient flow. By Assumption 
\ref{ass:stabilityPlant}, we have that $\underline{\lambda}(P_{\sigma})|\tilde{x}|^2\leq V_2(\tilde{x})\leq\overline{\lambda}(P_{\sigma})|\tilde{x}|^2$. Also, by Assumptions 
\ref{ass:lipschitzGradient} and \ref{ass:PLinequality}, we have that $\frac{\mu}{2}|\tilde{u}|^2\leq V_1(\tilde{u},\omega)\leq \frac{\ell}{2}|\tilde{u}|^2$, for all $\omega\in\Gamma$. Therefore, $\underline{a}|\tilde{z}_{\sigma}|^2\leq V_{\sigma}(\tilde{z},\omega)\leq \bar{a}_{\sigma}|z|^2$, for all $\omega\in \Gamma$, with 
$\bar a_{\sigma}$ and $\underline a_{\sigma}$ given by \eqref{barasigma} and \eqref{underbarasigma}.

\vspace{-.2cm}
\blue{
\textbf{Step 2.} Next, we show that for each fixed mode $\sigma\in\mathcal{S}$, and outside a neighborhood of the origin that is proportional in size to $\green{|\dot{\omega}|}$, the function $\green{V_\sigma}$ decreases along the trajectories of \eqref{errorplantdynamics} and \eqref{errorcontroldynamics} at an exponential rate. In particular, note that
}
\begin{align*}
\frac{1}{\eta}\nabla_{\tilde{u}}V_1^\top \dot{\tilde{u}}&=-\psi_c(0,\tilde{u},\omega)^\top \Big(\psi_c(\tilde{x},\tilde{u},\omega)-\frac{1}{\eta}\dot{\overbrace{u^*(\omega)}}\Big)\\
&\leq 
-\norm{\psi_c(0,\tilde{u},\omega)}^2 +a_2 \norm{\psi_c(0,\tilde{u},\omega)}  \norm{\tilde x}\\
&~~~+\frac{1}{\eta}\psi_c(0,\tilde{u},\omega)^\top \frac{\partial u^*}{\partial\omega}\dot{\omega},
\end{align*}
where the last inequality follows by Assumption \ref{ass:lipschitzGradient} with $a_2:=\ell_y  \norm{C} \norm{G}$.
%
Similarly, we have that

\vspace{-0.8cm}
\begin{small}
\begin{align}\label{dotV1w}
\nabla_{\omega} V_1^\top \dot{\omega}=\left(\nabla_{\omega} \phi(\tilde{u}+u^*(\omega),\omega)-\nabla_{\omega}\phi(u^*(\omega),\omega)\right) ^\top \dot{\omega}.
\end{align}
\end{small}
\vspace{-0.8cm}

Using the structure of $\phi$, and the chain rule:

\vspace{-0.8cm}
\begin{small}
\begin{align}\label{chainrule1}
\nabla_{\omega}\phi&=
\frac{\partial u^*}{\partial \omega}^\top \nabla\phi_u+H^\top \nabla\phi_y+\frac{\partial u^*}{\partial \omega}^\top G^\top \nabla\phi_y\notag\\
&=\frac{\partial u^*}{\partial \omega}^\top \nabla_u\phi+H^\top\nabla\phi_y.
\end{align}
\end{small}
\vspace{-0.8cm}

Since $\nabla_u\phi(u^*(\omega))=0$, and using the definition of $\psi_c(0,\tilde{u},\omega)$,  it follows that 

\vspace{-0.8cm}
\begin{small}
\begin{align*}
\nabla_{\omega} V_1^\top \dot{\omega} 
&\leq  -\left(\frac{\partial u^*}{\partial \omega}^\top\psi_c(0,\tilde{u},\omega)\right)^\top \dot{\omega}+\ell_y|H||G||\tilde{u}||\dot{\omega}|.
\end{align*}
\end{small}
\vspace{-0.8cm}

where we used again Assumption \ref{ass:lipschitzGradient}. Since $\dot{V}_1=\nabla_{\tilde{u}} V_1^\top \dot{\tilde{u}}+\nabla_{\omega} V_1^\top \dot{\omega}$, we have

\vspace{-0.8cm}
\begin{small}
\begin{align*}
\frac{1}{\eta}\dot{V}_1&\leq-\norm{\psi_c(0,\tilde{u},\omega)}^2+a_2 \norm{\psi_c(0,\tilde{u},\omega)}  \norm{\tilde x}\\
&~~~~+\frac{1}{\eta}\psi_c(0,\tilde{u},\omega)^\top \frac{\partial u^*}{\partial\omega}\dot{\omega}-\frac{1}{\eta}\left(\frac{\partial u^*}{\partial \omega}^\top\psi_c(0,\tilde{u},\omega)\right)^\top \dot{\omega}\\
&~~~~+\frac{\ell_y}{\eta}|H||G|\frac{|\psi_c(0,\tilde{u},\omega)|}{\mu}|\dot{\omega}|,\\
&\leq-\norm{\psi_c(0,\tilde{u},\omega)}^2+a_2 \norm{\psi_c(0,\tilde{u},\omega)}  \norm{\tilde x}\\
&~~~~+\frac{\ell_y}{\eta\mu}|H||G||\psi_c(0,\tilde{u},\omega)||\dot{\omega}|,
\end{align*}
\end{small}
\vspace{-0.8cm}

where the last term follows from the quadratic growth inequality. Using $\norm{\nabla_u \phi}^2 \geq 2 \mu V_1$, we obtain:

\vspace{-0.8cm}
\begin{small}
\begin{align*}
\frac{1}{\eta} \dot V_1 \leq - a_1
\norm{\psi_c(0,\tilde{u},\omega)}^2  +
a_2  \norm{\psi_c(0,\tilde{u},\omega)}  \norm{\tilde x}-a_3 V_1(\tilde{u},\omega),
\end{align*}
\end{small}
\vspace{-0.8cm}

where $a_1=(1-\kappa)$ and $a_3=\kappa\mu$,
which holds whenever $\norm{\psi_c(0,\tilde{u},\omega)} \geq \frac{2 \ell_y \norm{H}\norm{G}}{
\kappa  \eta \mu}\green{\norm{\dot \omega}}$. Since
$\norm{\psi_c(0,\tilde{u},\omega)} \geq \mu \norm{\tilde{u}}$, a sufficient condition for the above inequality to hold is
$\norm{\tilde{u}} \geq \frac{2 \ell_y \norm{H}\norm{G}}{
\kappa \eta \mu^2} \green{\norm{\dot \omega}}$.

On the other hand, for each fixed $\sigma\in\mathcal{S}$, we have:

\vspace{-0.8cm}
\begin{small}
\begin{align}
\frac{1}{\eta}
\dot V_2&=  
2 \tilde x^\tsp P_{\sigma} \left(\frac{1}{\eta} A_{\sigma} \tilde x + A_{\sigma}^\inv B_{\sigma} 
\psi_c(\tilde x, \tilde{u}, \omega) + \frac{1}{\eta} A_{\sigma}^\inv E_{\sigma} \dot \omega \right) \notag\\
&\leq - \frac{\underline \lambda(Q_{\sigma})}{\eta}  \norm{\tilde x}^2 + 
2 \norm{P_{\sigma} A_{\sigma}^\inv B_{\sigma}} \norm{\tilde x} \norm{\psi_c(\tilde x, \tilde{u}, \omega)} \notag\\
& \quad\quad+ \frac{2}{\eta} \norm{P_{\sigma} A_{\sigma}^\inv E_{\sigma}} \norm{\tilde x} \norm{ \dot \omega},\label{boundV2aux}
\end{align}
\end{small}
\vspace{-0.8cm}

where the inequality follows from Assumption \ref{ass:stabilityPlant}. Since $|\psi_c(\tilde x, \tilde{u},\omega)| \leq \norm{\psi_c(0, \tilde{u},\omega)} + a_2\norm{\tilde x}$, and  $V_2(\tilde x,\sigma) \leq \bar \lambda (P_{\sigma}) \norm{\tilde x}^2$, we can further upper bound \eqref{boundV2aux} as:
%
%
\begin{align*}
\frac{1}{\eta} \dot V_2 \leq 
-\frac{b_1}{\eta}\norm{\tilde x}^2+ b_2 \norm{\tilde x} 
\norm{\psi_c(0,\tilde{u},\omega)}+ b_3 \norm{\tilde x}^2 -b_4 V_2(\tilde x,\sigma)
\end{align*}
which holds whenever $\norm{\tilde x} \geq 
\frac{4 \norm{P_{\sigma} A_{\sigma}^\inv E_{\sigma}}}{\kappa \underline \lambda(Q_{\sigma})} 
\green{\norm{\dot \omega}}$, where $b_1:=(1-\kappa) \underline \lambda(Q_{\sigma})$, $b_2:=2 \norm{P_{\sigma}A_{\sigma}^\inv B_{\sigma}}$, $b_3:=2 a_2 \norm{P_{\sigma}A_{\sigma}^\inv B_{\sigma}}$, and $b_4:=\frac{\kappa}{2\eta} 
\frac{\underline \lambda(Q_{\sigma})}{\bar \lambda(P_{\sigma})}$.

\textbf{Step 3.} \blue{ Next, using the bounds on $\dot{V}_1$ and $\dot{V}_2$, we obtain that for each fixed $\sigma\in\mathcal{S}$ the function $V_{\sigma}$ satisfies the following bound outside a $|\dot{\omega}|$-neighborhood of the origin: }

\vspace{-1.0cm}
\begin{small}
\begin{align}
\label{eq:Vdot_Gradient}
\dot V_{\sigma}&\leq
- \xi^\tsp \Lambda_{\sigma} \xi 
- \green{\frac{\kappa}{2} \min \left\{2 \eta \mu, 
\frac{\underline \lambda (Q_{\sigma})}{\bar \lambda (P_{\sigma})} \right\}} V_{\sigma}(\tilde{z},\omega),
\end{align}
\end{small}
\vspace{-0.7cm}

where $\xi = (\norm{\nabla \phi(\tilde{u}+u^*(\omega),\omega)}, \norm{\tilde x})$
and $\Lambda_{\sigma}$ is a $2\times 2$ symmetric matrix with entries
$\Lambda_{11}=(1-\theta) a_1$,
$\Lambda_{12}=\Lambda_{12}=-\frac{1}{2}((1-\theta) a_2 + \theta b_2)$, 
and $\Lambda_{22}=\theta (b_1/ \eta - b_3)$.
Since the largest $\theta\in(0,1)$ that guarantees $\Lambda \succ 0$
is as in \eqref{eq:aBaraUnderbarGradient}, we conclude that $\Lambda\succ0$ when $\eta_{\sigma} \in(0, \bar \eta_{\sigma})$, with $\bar{\eta}_{\sigma}$ as in \eqref{eq:etaBoundGradient}.

\vspace{-.1cm}
\textbf{Step 4.} 
\blue{
Finally, we incorporate the switches of the plant which are governed by the states $(\sigma,\tau)$ generated by \eqref{eq:dwellTimeDynamics_hybrid}.} Using $V_{\sigma}$, we consider the extended function $W(\tilde{\vartheta}) = e^{\varrho \tau} V_{\sigma}(\tilde{z},\omega)$, where $\varrho>0$ and $\tilde{\vartheta}:=(\tilde{z},\sigma,\tau,\omega)$. We will show the existence of $\varrho$ such that $W$ is a hybrid ISS Lyapunov function for the HDS with dynamics \eqref{eq:dwellTimeDynamics_hybrid}, \eqref{exosystemdynamics}, and \eqref{errorplantdynamics}-\eqref{errorcontroldynamics}, with respect to the ``input'' $|\dot{\omega}|$, and the compact set $\mathcal{A}_1:=\left\{\tilde{\vartheta}:\tilde{z}=0,~\sigma\in\mathcal{S},~\tau\in[0,N_0],~\omega\in\Gamma\right\}$. Indeed, note that for all $\tilde{z}\in\mathbb{R}^{n+m}$, all $(\tau,\sigma)\in[0,N_0]\times\mathcal{S}$, and $\omega\in\Gamma$, we have that
%
$\min_{\sigma} \underline a_{\sigma} \norm{\tilde{z}}^2 
\leq 
W(\tilde{\vartheta})\leq e^{\varrho N_0}\max_{\sigma} \bar a_{\sigma}  \norm{\tilde{z}}^2$,
where $\underline a_{\sigma}, \bar a_{\sigma}$ are as in 
\eqref{eq:aBaraUnderbarGradient}. During flows, we have that:
\begin{align*}
\hspace{-.35cm} \dot W 
= 
e^{\varrho \tau} (\varrho V_{\sigma} \dot \tau +  \dot V_{\sigma})  \leq 
\left(\frac{\varrho}{\tau_d} - \min_{\sigma} b_{\sigma} \right) W(\tilde{\vartheta}),
\end{align*}
where $b_{\sigma}$ is as in \eqref{bsigma}, which holds whenever
$\norm{\tilde{u}} \geq  \green{\frac{2 \ell_y \norm{H}\norm{G}}{
\kappa \eta \mu^2} \norm{\dot \omega}}$ and
$\norm{\tilde x} \geq 
\green{\frac{4 \norm{P_\sigma A_\sigma^\inv E_\sigma}}{\kappa \underline \lambda(Q_\sigma)} \norm{\dot \omega}}$.
Hence, $W$ decreases during flows if $\varrho <  \tau_d \min_{\sigma} b_{\sigma}$. Next, note that at each jump we have that $\tilde{x}^+=\tilde{x}$, $\tilde{u}^+=\tilde{u}$, $\sigma^+\in\mathcal{S}$, $\omega^+=\omega$, and $\tau^+=\tau-1$. Hence, $W(\tilde{\vartheta}^+) \leq e^{\varrho \tau^+} \max_{\sigma} V_{\sigma}(\tilde{z}^+,\omega^+)$ and $W(\tilde{\vartheta}^+) \leq e^{-\varrho + \ln(\max_{\sigma}\bar{a}_{\sigma}) - \ln(\min_{\sigma}\underline{a}_{\sigma})} W(\tilde{\vartheta})$. It follows that if 
$\varrho > \ln\left(\max_{\sigma}\bar a_{\sigma} /\min_{\sigma} \underline a_{\sigma}\right)$, the function $W$ also decreases during jumps. It follows that if $\tau_d > \frac{ \ln (\max_{\sigma} \bar a_{\sigma} / \min_{\sigma}\underline a_{\sigma})}{\min_{\sigma} b_{\sigma}}$, then $W$ is a hybrid ISS Lyapunov function for the closed-loop system with ``input'' $\dot{\omega}$. The quadratic bounds and the periodic nature of the jumps imply exponential input-to-state stability of $\mathcal{A}_1$, with input $|\dot{\omega}|$.  \QEDBL

\subsection{Proof of Theorem \ref{thm:EISS-Nesterov}}
\label{sec:proofEISS_Nesterov}
We follow similar steps to the proof of Theorem \ref{thm:EISS-Gradient}, now incorporating the jumps of the hybrid controller.

\noindent 
\textbf{Step 1.} 
\blue{
Let $\tilde{u}_1:=u_1-u^*(\omega)$, $\tilde{u}_2:=u_2-u^*(\omega)$, $\tilde{u}_3=u_3$, and $\tilde{x}:=x+A_{\sigma}^\inv B_{\sigma} (\tilde{u}_1+u^*(\omega)) + A_{\sigma}^\inv E_{\sigma} \omega$. The error dynamics of $\tilde{u}=(\tilde{u}_1,\tilde{u}_2,\tilde{u}_3)$ are given by $\dot{\tilde{u}}_3=\frac{\eta}{2}$, and
}
\begin{align}
\dot{\tilde{u}}_1&=\dfrac{2\eta}{\tilde{u}_3} (\tilde{u}_2-\tilde{u}_1)-\dot{\overbrace{u^*(\omega)}},\notag\\
\dot{\tilde{u}}_2
&=2k\eta \tilde{u}_3  \psi_c(\tilde{x},\tilde{u}_1,\omega)-\dot{\overbrace{u^*(\omega)}},\label{errorcontroldynamics3}
\end{align}
where $\psi_c$ was defined in \eqref{errorplantdynamics}. Also, we have
\begin{align}\label{errorplantdynamics3}
\dot{\tilde{x}}=A_{\sigma}\tilde{x}+\frac{2\eta}{\tilde{u}_3}A_{\sigma}^{-1}B_{\sigma}(\tilde{u}_2-\tilde{u}_1)+A_{\sigma}^{-1}E_{\sigma}\dot{\omega}.    
\end{align}
We consider the Lyapunov-like function:
\begin{align}
\label{eq:LyapunovFunction_Nesterov}
V_{\sigma}(\tilde{z},\omega) = \frac{(1-\green{\theta_\sigma})}{\green{\eta_\sigma}}  V_1(\tilde{u},\omega) + \frac{\green{\theta_\sigma}}{\green{\eta_\sigma}} V_2(\tilde x,\sigma),~
\end{align}
where $\theta_{\sigma}\in(0,1)$, $\tilde{z}:= (\tilde{x}, \tilde{u})$, $\tilde{u}:=(u_1, u_2, u_3)$, 
$V_1(\tilde{u},\omega)=\frac{1}{4}|\tilde{u}_2-\tilde{u}_1|^2+\frac{1}{4}|\tilde{u}_2|^2+k\tilde{u}_3^2(\phi(\tilde{u}_1+u^*(\omega),\omega)-\phi(u^*(\omega),\omega))$
and $V_2(\tilde x,\sigma) = \tilde x^\tsp P_{\sigma} \tilde x$. Note that $V_1$ is a Lyapunov function for the accelerated hybrid gradient controller acting on \emph{static} maps \cite{JIP-NL:19}. The individual components of $V_{\sigma}$ satisfy $\underline \lambda(P_{\sigma}) \norm{\tilde x} \leq V_2(\tilde x,\sigma) \leq \bar \lambda (P_{\sigma})\norm{\tilde x}$ with $\bar a_{\sigma}, \underline a_{\sigma}$ as in \eqref{eq:aBaraUnderbarGradient}; and $\green{\frac{1+2k \ell \delta^2}{4}} 
\norm{\tilde{u}}_{\mathcal{A}_{\tilde{u}}}^2 
\leq 
V_1(\tilde{u},\omega) 
\leq  
\green{\frac{1 + 2 k  \ell \Delta^2}{2}} 
\norm{\tilde{u}}_{\mathcal{A}_{\tilde{u}}}^2$, for all $\omega\in\Gamma$, where $\mathcal{A}_{\tilde{u}}=\{0\}\times\{0\}\times[\delta,\Delta]$, and where we used Assumptions \ref{ass:lipschitzGradient} and \ref{ass:strongConvexity}. It follows that for each $\sigma\in\mathcal{S}$, $V_{\sigma}$ satisfies $\underline{a}_{\sigma}|\tilde{z}|^2\leq V_{\sigma}(\tilde{z})\leq \bar{a}_{\sigma}|\tilde{z}|^2$, for all $\tilde{z}\in\mathbb{R}^{n+m}$ and all $\omega\in\Gamma$, where $\tilde{z}:=(\tilde{x},\tilde{u})$.
Next, we show that $V_{\sigma}$ decreases during the flows of $\tilde{u}$ and $\tilde{x}$. In particular, since $\psi_c(0,\tilde{u}_1,\omega)=-\nabla_{u_1} \phi(\tilde{u}_1+u^*(\omega),\omega)$, we have
\vspace{-0.5cm}
\begin{small}
\begin{align*}
&\frac{1}{\eta}\nabla_{\tilde{u}_1}V_1^\top \dot{\tilde{u}}_1=
-\frac{1}{\tilde{u}_3} \norm{\tilde{u}_2-\tilde{u}_1}^2 +\frac{1}{2\eta}(\tilde{u}_2-\tilde{u}_1)^\top\dot{\overbrace{u^*(\omega)}}\\
& ~~~-2k\tilde{u}_3\psi_c(0,\tilde{u}_1,\omega)^\top (\tilde{u}_2-\tilde{u}_1) 
+\frac{k}{\eta}\tilde{u}_3^2\psi_c(0,\tilde{u}_1,\omega)^\top\dot{\overbrace{u^*(\omega)}},
\end{align*}
\end{small}
\vspace{-0.5cm}
and
\vspace{-0.3cm}
\begin{small}
\begin{align*}
\frac{1}{\eta}\nabla_{u_2} V_1^\top  \dot {\tilde{u}}_2 
&=\frac{1}{2\eta}(2\tilde{u}-\tilde{u}_1)^\top\dot{\tilde{u}}_2\\
& \leq k \tilde{u}_3  (2 \tilde{u}_2 - \tilde{u}_1)^\tsp \psi_c(0,\tilde{u}_1,\omega)+ 2 \bar a_2 \norm{\tilde x}|\tilde{u}_1-\tilde{u}_2|\\
&~~~+\bar a_2 |\tilde{x}||\tilde{u}_1|-\frac{(2\tilde{u}_2-\tilde{u}_1)^\top}{2\eta} \dot{\overbrace{u^*(\omega)}},
\end{align*} 
\end{small}\noindent 
where $\bar a_2:= k\Delta\ell_y \norm{C}\norm{G}$. Also, using \eqref{dotV1w} and \eqref{chainrule1}:
\begin{align*}
\frac{1}{\eta}\nabla_{\omega} V_1^\top \dot{\omega}
&\leq -\frac{k\tilde{u}_3^2}{\eta}\left(\frac{\partial u^*}{\partial \omega}^\top\psi_c(0,\tilde{u}_1,\omega)\right)^\top \dot{\omega}\\
&~~~~+\frac{k\tilde{u}_3^2}{\eta}\ell_y|H||G||\tilde{u}_1||\dot{\omega}|\leq \green{\frac{\bar a_3}{\eta}} |\tilde{u}|_{\mathcal{A}_{\tilde{u}}}|\dot{\omega}|,
\end{align*}
with  \green{$\bar a_3:= k \Delta^2 (\ell|\frac{\partial u^*}{\partial\omega}|+\ell_y|H||G|)$}, where we used Lipschitz continuity of $\nabla_{u_1}\phi(\cdot)$. Combining the above bounds and using $\frac{1}{\eta}\nabla_{\tilde{u}_3}V_1^\top\dot{\tilde{u}}_3=k \tilde{u}_3(\phi(\tilde{u}_1+u^*(\omega),\omega)-\phi(u^*(\omega),\omega))$, we obtain 
\vspace{-.5cm}
\begin{small}
\begin{align}
&\frac{1}{\eta}\nabla V_1^\top  \dot{\tilde{u}}\leq 
-\tilde{u}_3^{-1}|\tilde{u}_2-\tilde{u}_1|^2+\bar a_2|\tilde{x}||\tilde{u}_1|+\green{\frac{1}{\eta}}\bar a_3|\dot{\omega}||\tilde{u}|_{\mathcal{A}_{\tilde{u}}} \notag\\
&~~~-k\tilde{u}_3\Big(\phi(u^*(\omega),\omega)  -\phi(\tilde{u}_1+u^*(\omega),\omega)
-\psi_c(0,\tilde{u}_1,\omega)^\top\tilde{u}_1\Big)
\notag\\
&~~~+2a_2 \norm{\tilde x}|\tilde{u}_1-\tilde{u}_2|
+\frac{1}{2\eta}\dot{\overbrace{u^*(\omega)}}^\top \left(-\tilde{u}_2+2k\tilde{u}^2_3\psi_c(0,\tilde{u}_1,\omega)\right).\label{boundV1thm33}
\end{align}
\end{small}
\vspace{-.6cm}
Using strong convexity to bound the fourth term,
%
%
%

\vspace{-0.3cm}
\begin{small}
\begin{align}
\label{eq:V1dot_nesterov_final}
\frac{1}{\eta}\nabla V_1^\top \dot{\tilde{u}}
&\leq - \min\{\frac{1}{\Delta}, \frac{k \delta \mu}{4}\}
\norm{\tilde{u}}_{\mathcal{A}_{\tilde{u}}}^2
+ \green{a_2} \norm{\tilde x} \vert \tilde{u} \vert_{\mathcal{A}_{\tilde{u}}}
+ \green{\frac{1}{\eta} a_3} \norm{\dot \omega} \vert \tilde{u} \vert_{\mathcal{A}_{\tilde{u}}} \nonumber\\
&\hspace{-1.1cm} \leq - \green{a_1}
\norm{\tilde{u}}_{\mathcal{A}_{\tilde{u}}}^2
+ \green{a_2} \norm{\tilde x} \vert \tilde{u} \vert_{\mathcal{A}_{\tilde{u}}}
\green{- \frac{\kappa}{2} \frac{2 \min \left\{\frac{1}{\Delta}, 
\frac{k \delta \mu}{4}\right\}}{1 + 2 k \Delta^2 \ell} V_1(\tilde u, \omega)} \hspace{-.1cm} 
\end{align}
\end{small}

\vspace{-0.7cm}\noindent
where $a_2:=4\sqrt{2}\bar a_2$, 
$a_1 := (1-\kappa) \green{\min \left\{\frac{1}{\Delta}, 
\frac{k \delta \mu}{4}\right\}}$,
$a_3:=\bar a_3+\frac{1}{2}\max_{\omega\in\Gamma}\left|\frac{\partial u^*}{\partial \omega}\right|\sqrt{2}\max\{1,2k\ell\Delta^2\}$, and where we used 
the inequalities 
$\vert \tilde{u} \vert_{\mathcal{A}_{\tilde{u}}}^2 \leq \norm{\tilde{u}_1}^2 + 2 \norm{\tilde{u}_2-\tilde{u}_1}^2 + 2 \norm{\tilde{u}_1}^2$ and
\green{$V_1 (\tilde u, \omega)\leq \frac{1}{2} (1 + k\ell\Delta^2) \vert \tilde{u} \vert^2_{\mathcal{A}_{\tilde{u}}}$},
and where the last inequality holds whenever
\green{$\vert \tilde{u} \vert_{\mathcal{A}_{\tilde{u}}} \geq \frac{2 a_3}{\kappa \eta \min \left\{\frac{1}{\Delta}, 
\frac{k \delta \mu}{4}\right\}} \norm{\dot \omega}$.}

\blue{
\textbf{Step 2.} Next, we show that the reset policy $r_0=1$ guarantees the decrease of $V_1$ during the discrete-time updates of the controller. In particular, in this case we have that $\Delta V_1:=V_1(\tilde{u}^+,\omega) - V_1(\tilde{u},\omega)$ satisfies
}

\vspace{-.5cm}
\begin{small}
\begin{align*}
 \Delta V_1 &= 
\frac{1}{4}\norm{\tilde{u}_1}^2 + 
k \delta^2(\phi(\tilde{u}_1+u(\omega)^*,\omega)-\phi(u^*(\omega),\omega)) -\frac{1}{4}\norm{\tilde{u}_2}^2 \nonumber\\
&\quad-\frac{1}{4}\norm{\tilde{u}_2 - \tilde{u}_1}^2 
- k\Delta^2(\phi(\tilde{u}_1+u^*(\omega),\omega)-\phi(u^*(\omega),\omega)) \nonumber\\
%
&\leq - c_s V_1(\tilde{u},\omega),
\end{align*}
\end{small}
\vspace{-.8cm}

where $c_s:=1-\frac{\delta^2}{\Delta^2} - \frac{1}{
2\mu \Delta^2 k}$, and where we used 
$\norm{\tilde{u}_1}^2 \leq \frac{2}{\mu} (\phi(\tilde{u}_1+u^*(\omega),\omega)-\phi(u^*(\omega),\omega))$ for the first inequality, and $\Delta^2 - \delta^2 \geq \frac{1}{2k\mu}$ for the second inequality. This establishes that $V_1(\tilde{u}^+,\omega) \leq (1-c_s) V_1(\tilde{u},\omega)$, with $c_s\in(0,1)$.

\blue{
\textbf{Step 3.} We now consider the function $V_2$ in~\eqref{eq:LyapunovFunction_Nesterov}, and for each fixed $\sigma$ we bound its 
evolution during flows and jumps triggered by the controller. In this case, we have
}

\vspace{-.8cm}
\begin{small}
\begin{align}
\frac{1}{\eta}
\nabla V_2^\top \dot{\tilde{x}}&= 
\frac{1}{\eta} \tilde x^\tsp (A_{\sigma}^\tsp P_{\sigma} + P_{\sigma}A_{\sigma}) \tilde x 
+ 2 \tilde x^\tsp P_{\sigma}A_{\sigma}^\inv B_{\sigma} \frac{(\tilde{u}_2-\tilde{u}_1)}{\tilde{u}_3}\notag\\
& \quad\quad + \frac{1}{\eta} 2 \tilde x^\tsp P_{\sigma}A_{\sigma}^\inv E_{\sigma} \dot \omega\notag\\
& \hspace{-1.1cm} \leq 
- \frac{1}{\eta} \underline \lambda(Q_{\sigma}) \norm{\tilde x}^2 
+ b_2 \norm{\tilde x} \norm{\tilde{u}}_{\mathcal{A}_{\tilde{u}}}+ 2 \eta^\inv \norm{P_{\sigma}A_{\sigma}^\inv E_{\sigma}} \norm{\tilde x} \norm{\dot \omega}\label{V2dotThm33} \hspace{-.2cm}
\end{align}
\end{small}
\vspace{-.9cm}

where \green{$b_2:=4\sqrt{2} \delta^\inv \norm{P_{\sigma}A_{\sigma}^\inv B_{\sigma}}$}. Using
$V_2(\tilde x,\sigma) \leq \bar \lambda(P_{\sigma}) \norm{\tilde x}^2$,
we obtain
\begin{align}\label{decreaseV2Thm3}
\hspace{-.3cm} \frac{1}{\eta}
\dot V_2 \leq
- \frac{1}{\eta} b_1
\norm{\tilde x}^2  
+ b_2 \norm{\tilde x} \norm{\tilde{u}}_{\tilde{\mc A}}- \frac{\kappa}{2\eta} \frac{\underline \lambda (Q_\sigma)}{\bar \lambda(P_{\sigma})} 
V_2(\tilde x, \sigma)
\end{align}
where $b_1:=(1-\kappa) \underline \lambda(Q_{\sigma})$, which holds only if $\norm{\tilde x} \geq \frac{4 \norm{P_{\sigma}A_{\sigma}^\inv E_{\sigma}}}{\kappa \underline \lambda(Q_{\sigma})} 
\sup_\tau \norm{\dot \omega_\tau}$. Since the state of the plant and its mode do not change during discrete-time updates of the controller, we have that $V_2(\tilde{x}^+,\sigma^+)=V_2(\tilde{x},\sigma)$.
%

\blue{
\textbf{Step 4.} Combining the estimates \eqref{eq:V1dot_nesterov_final} and \eqref{decreaseV2Thm3}, we obtain that for each fixed mode $\sigma$, and outside a $|\dot{\omega}|$-neighborhood of $\{0\}\times\mathcal{A}_{\tilde{u}}$, $V_{\sigma}$ satisfies $\dot V_{\sigma}(\tilde{z}) \leq
- \xi^\tsp \Lambda_{\sigma} \xi - 
b_\sigma
V_{\sigma}(\tilde{z})$, where $\xi = (\vert \tilde{u} \vert_{\mathcal{A}_{\tilde{u}}}, \norm{\tilde x})$, $b_\sigma$ is as in \eqref{eq:coefficientsNesterov-d},
and $\Lambda_{\sigma}$ is a $2\times 2$ symmetric matrix with entries
$\Lambda_{11}=(1-\theta) a_1$,
$\Lambda_{12}=\Lambda_{12}=-\frac{1}{2}((1-\theta) \tilde{a}_2 + \theta b_2)$, and
$\Lambda_{22}=\theta b_1/\eta$. It follows that when $\eta < \hat{\eta}_{H}$, with $\hat{\eta}_{H}$ given by \eqref{eq:etaBoundNesterov}, the matrix $\Lambda_{\sigma}$ is positive definite. Similarly, during jumps triggered by the controller: }
\begin{align}\label{eq:VjumpsGradient}
V_{\sigma^+}(\tilde{z}^+) \leq \min \{1-c_s, 1\} V_{\sigma}(\tilde{z}) \leq V(\tilde{z}),
\end{align} 
which implies that $V_{\sigma}$ does not increase.

\vspace{-.1cm}
\textbf{Step 6.} Finally, we incorporate the switches into the error plant dynamics \eqref{errorplantdynamics}, where $(\sigma,\tau)$ are generated by \eqref{eq:dwellTimeDynamics_hybrid}. Using $V_{\sigma}$, we consider the extended function $W(\tilde{\vartheta}) = e^{\varrho \tau} V_{\sigma}(\tilde{z})$, where $\varrho>0$ and $\tilde{\vartheta}:=(\tilde{z},\sigma,\tau,\omega)$, with $\tilde{z}=(\tilde{x},\tilde{u}_1,\tilde{u}_2,\tilde{u}_3)$. We will show the existence of $\varrho$ such that $W$ is a hybrid ISS Lyapunov function for the HDS with error dynamics \eqref{errorcontroldynamics3}-\eqref{errorplantdynamics3}, exosystem \eqref{exosystemdynamics}, and switching generator \eqref{eq:dwellTimeDynamics_hybrid}, with respect to the compact set:
\begin{equation*}
\mathcal{A}_3:=\left\{\tilde{\vartheta}:\tilde{x}=0,\tilde{u}\in\mathcal{A}_{\tilde{u}},~\sigma\in\mathcal{S},~\tau\in[0,N_0],~\omega\in\Gamma\right\}   
\end{equation*}

and with ``input'' $\dot{\omega}$. Indeed, note that $\min_{\sigma} \underline a_{\sigma} \norm{\tilde{\vartheta}}^2 
\leq W(\tilde{\vartheta})\leq e^{\varrho N_0}\max_{\sigma} \bar a_{\sigma}  \norm{\tilde{\vartheta}}^2$, where $\underline a_{\sigma}, \bar a_{\sigma}$ are as in Theorem \ref{thm:EISS-Nesterov}. During flows of the system, we have that:

\vspace{-0.8cm}
\begin{small}
\begin{align*}
\hspace{-.35cm} \nabla W^\top\dot{\tilde{\vartheta}} 
= e^{\varrho \tau} (\varrho V_{\sigma}(\tilde{z}) \dot \tau +  \dot V_{\sigma}(\tilde{z}))  
\leq 
\left(\frac{\varrho}{\tau_d} - \min_{\sigma} b_{\sigma} \right) W(\tilde{\vartheta}),
\end{align*}
\end{small}
\vspace{-0.7cm}

which holds for all $\vert u \vert_{\mathcal{A}_{\tilde{u}}} \geq \frac{2 \bar{c}_1}{\kappa \frac{1}{3} \min \{ \tilde{u}_3^\inv, \frac{\tilde{u}_3 \mu}{2}\}} 
\norm{\dot \omega}$ and 
$\norm{\tilde x} \geq \frac{4 \norm{PA^\inv E}}{\kappa \underline \lambda(Q)} \norm{\dot \omega}$. Similarly, during jumps of the form $\tilde{\varrho}  \in D_1$, we have that $\tilde{x}^+=\tilde{x}$, 
$\tilde{u}^+=\tilde{u}$, $\sigma^+\in\mathcal{S}$, and $\tau^+=\tau-1$, and thus $W(\tilde{\varphi}^+) \leq e^{-\varrho + \ln(\max_{\sigma}\bar{a}_{\sigma}) - \ln(\min_{\sigma}\underline{a}_{\sigma})} W(\tilde{\varphi})$.
%
%
%
During jumps of the form $\tilde{\varrho} \in D_2$, we have
%
$W(\tilde{\varphi}^+) = e^{\varrho \tau} V_{\sigma}(\tilde{z}^+)  \leq W(\tilde{\varphi})$,
where we used \eqref{eq:VjumpsGradient}. Therefore, to guarantee that $W$ does not increase during jumps it suffices to have $\varrho >   \ln(\max_{\sigma}\bar a_{\sigma} / \min_{\sigma}\underline a_{\sigma})$ and $\varrho <  \min_{\sigma} b_{\sigma} \tau_d$. Combining the upper and lower inequalities on $\varrho$ we conclude that we need $\tau_d > \frac{\ln ( \max_{\sigma}\bar a_s /\min_{\sigma} \underline a_{\sigma})}{\min_{\sigma} b_{\sigma}}$, which establishes the result.  \QEDBL

\vspace{-.2cm}
\subsection{Proof of Theorem \ref{thm:EISS-Nesterov-Practical}}
\vspace{-.35cm}
We follow similar steps as in the proof of Theorem \ref{thm:EISS-Nesterov}. Since $\omega$ is now constant, we have that $\dot{\omega}=0$, and since we focus on a single mode of the plant we drop the subscript $\sigma$. The error dynamics \eqref{errorcontroldynamics3} become $\dot{\tilde{u}}_3=\frac{\eta}{2}$ and
\begin{align}\label{errordynamicsthm32}
\dot{\tilde{u}}_1=\dfrac{2\eta}{\tilde{u}_3} (\tilde{u}_2-\tilde{u}_1),~\dot{\tilde{u}}_2=2k\eta \tilde{u}_3  \psi_c(\tilde{x},\tilde{u}_1,\omega),
\end{align}
while the dynamics of $\tilde{x}$ are still given by \eqref{errorplantdynamics3}. By construction, and under Assumption \ref{ass:reverseLipschitz}, the function $V$ defined in \eqref{eq:LyapunovFunction_Nesterov} is radially unbounded and positive definite with respect to the compact set $\{0\}\times\mathcal{A}_{\tilde{u}}$. Along the dynamics \eqref{errordynamicsthm32}, the function $V_1$ still satisfies  \eqref{boundV1thm33}. Using convexity and the Lipschitz property of $\nabla_u\phi$, we obtain:
\begin{align*}
\frac{1}{\eta}
\nabla V_1^\top\dot{\tilde{u}} & \leq 
- 2 \Delta^\inv \norm{\tilde{u}_2-\tilde{u}_1}^2  
+ a_2 \norm{\tilde x} \norm{\tilde{u}_1}+2a_2 \norm{\tilde x}|\tilde{u}_1-\tilde{u}_2|\\
&~~~~-a_4\norm{\nabla \phi(\tilde{u}_1+u^*)}^2,
\end{align*}
where $a_4:=\frac{\delta k}{2\ell}$. Similarly, from \eqref{V2dotThm33}, we now have that
\begin{align}
\label{eq:V2dot_practical_final} 
\frac{1}{\eta}
\nabla V_2^\top \dot{\tilde{x}}
&\leq 
-\frac{1}{\eta} \underline \lambda (Q) \norm{\tilde x}^2 
+ b_2 \norm{\tilde x}\norm{\tilde{u}}_{\mathcal{A}_{\tilde{u}}}.
\end{align}
where $b_2:=2 \delta^{-1} \norm{PA^\inv B}$. It follows that $V$ satisfies
\begin{align}
\nabla V^\top\dot{\tilde{u}} \leq & - \xi^\tsp  \Lambda  \xi -a_4 (1-\theta) \norm{\nabla \phi(\tilde{u}_1+u^*)}^2\notag\\  
&- \kappa \theta \frac{\underline \lambda(Q)}{\eta} \norm{\tilde x}^2 + (1-\theta)a_2 \norm{\tilde x}  \norm{\tilde{u}_1},
\end{align}
with $\xi = ( \norm{\tilde{u}_2-\tilde{u}_1}, \norm{\tilde x})$, and where $\Lambda \in \real^{2 \times 2}$ is a symmetric matrix 
with entries given by $\Lambda_{11} = (1-\theta) \Delta^\inv$,
$\Lambda_{12}=\Lambda_{21}= -\frac{1}{2}((1-\theta)2a_2 + \theta b_2)$, and
$\Lambda_{22} = - \theta (1-\kappa) \underline \lambda(Q)/2$. We now consider two possible scenarios:

\vspace{-.35cm}
\textsl{Case 1:~} Suppose that
$- \kappa \theta \frac{\underline \lambda(Q)}{\eta} \norm{\tilde x}^2 
+ (1-\theta) a_2 \norm{\tilde x}  \norm{\tilde{u}_1} \leq 0$.
In this case, we have
$\dot V \leq - \xi^\tsp  \Lambda  \xi  - (1-\theta) a_4 \norm{\nabla \phi(\tilde{u}_1+u^*)}^2$; since $\Lambda \succ 0$, and using the convexity of $\phi$, we conclude that $\dot V <0$ for all $\eta \in (0, \bar \eta_s)$ and $\tilde{z}\neq (0,u^*)$.

\vspace{-.35cm}
\textit{Case 2:~} Suppose 
$- \kappa \theta \frac{\underline \lambda(Q)}{\eta} \norm{\tilde x}^2 
+ (1-\theta) a_3 \norm{\tilde x}  \norm{\tilde{u}_1} > 0$.
In this case, we note that $\eta \in (0, \bar \eta_s)$ implies 
$\sqrt{\frac{1-\theta}{\theta} \frac{4 \ell a_3^2 \eta}{\kappa^2 b_1 \delta k}} \leq \ell_0$ and, by combining this observation with 
Assumption 
\ref{ass:reverseLipschitz} we conclude that $\norm{\nabla \phi(\tilde{u}_1+u^*)} \geq \sqrt{\frac{1-\theta}{\theta} \frac{4 \ell a_3^2 \eta}{\kappa^2 b_1 \delta k}} \norm{\tilde{u}_1}$ and thus $\dot V\leq - \xi^\tsp \Lambda \xi - \kappa \theta \frac{\underline \lambda(Q)}{\eta} \norm{\tilde x}^2 - (1-\kappa) (1-\theta) \frac{k \delta }{4 \ell} \norm{\nabla \phi(\tilde{u}_1+u^*)}^2$. Since $\Lambda \succ 0$ when $\eta \in (0, \bar \eta_s)$, we have that  $\dot V<0$ outside a neighborhood of $\{0\}\times \mathcal{A}_{\tilde{u}}$. Finally, we show that if $\Delta>\delta$, the Lyapunov function $V$ does not increase during jumps of controller. Indeed, the reset policy $r_0=0$ implies that $V(\tilde{z}^+) - V(\tilde{z})= \frac{k}{4} (\delta^2 - \Delta^2)
\left(\phi(\tilde{u}_1+u^*)-\phi(u^*)\right)\leq0$, where the inequality holds because $\delta < \Delta$. The strong decrease of $V$ during flows outside of a neighborhood of $\{0\}\times\mathcal{A}_{\tilde{u}}$, the non-increase during jumps, and the periodic hybrid time domain of the solutions guarantee uniform convergence of $(\tilde{x},u_1,u_2)$ from compact sets to a neighborhood of $(0,u^*,u^*)$ via \cite[Ch.8]{RG-RGS-ART:12}, which establishes item a) of Theorem 3.2. Item b) follows by the fact that $\dot{V}\leq0$ outside a neighborhood of $\{0\}\times\mathcal{A}_{\tilde{u}}$, which implies that in this set $V(t,j)\leq V(\underline{t}_j,j)$ for each $(t,j)$ in the domain of the solution, and the fact that by construction of \eqref{eq:LyapunovFunction_Nesterov} we have $ku_3(t,j)^2(\phi(u_1(t,j))-\phi^*)\leq V(t,j)$, which leads to the bound of the theorem with $c_j:=k^{-1}V(\underline{t}_j,j)$.  \QEDBL

\vspace{-.3cm}

\section{Numerical Examples}
\label{sec:6}

\vspace{-.3cm}

\begin{figure}[tb]
\centering
\subfigure[]{\includegraphics[width=.5\columnwidth]{%
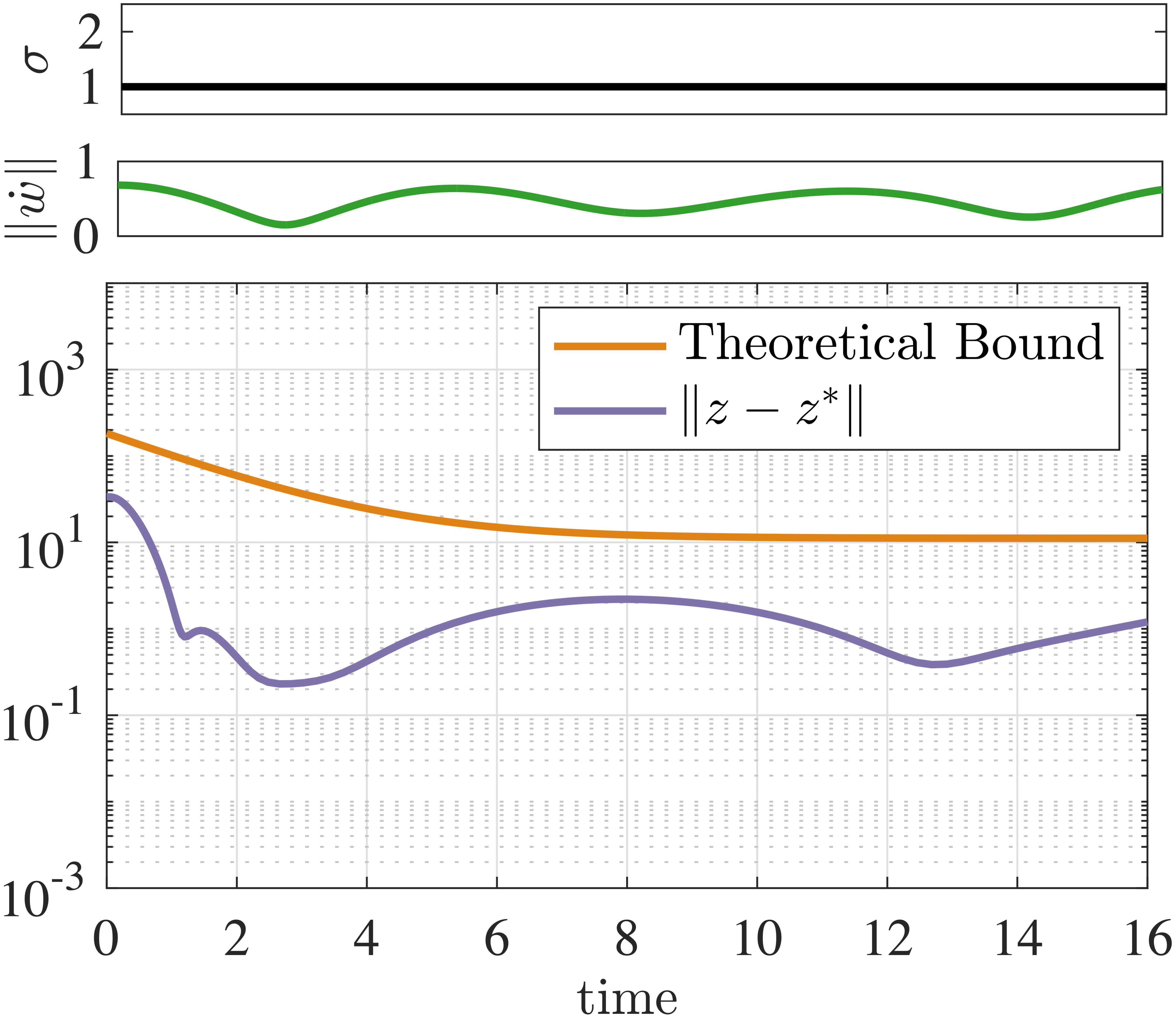}}%
\subfigure[]{\includegraphics[width=.5\columnwidth]{%
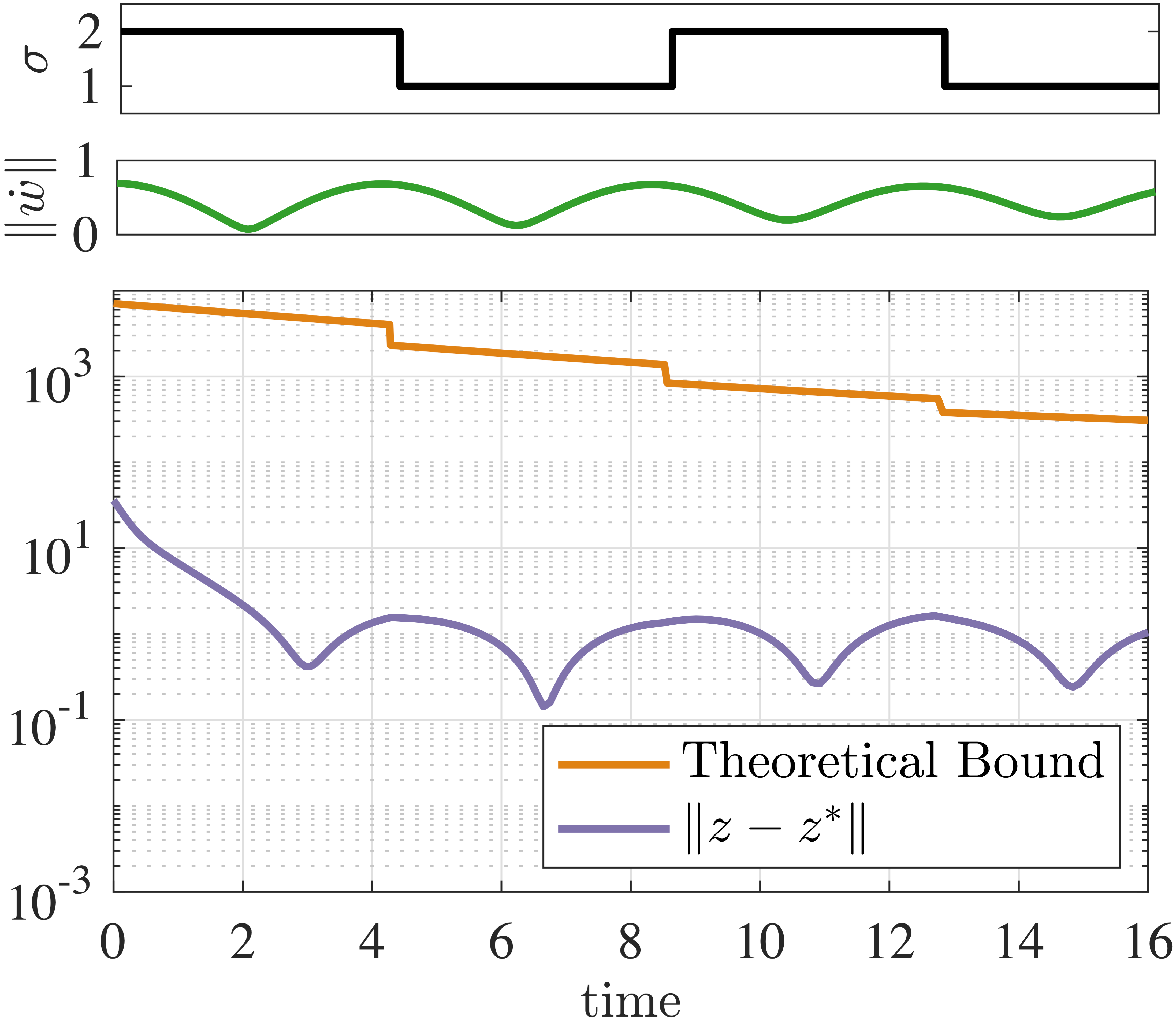}}
\caption{Tracking error of gradient flow controller 
\eqref{eq:gradient_controller}. (a) No plant switching. (b) Plant switches 
between two modes.}
\label{fig:thm-3-3} 
\end{figure}

\begin{figure}[tb]
\includegraphics[width=\columnwidth]{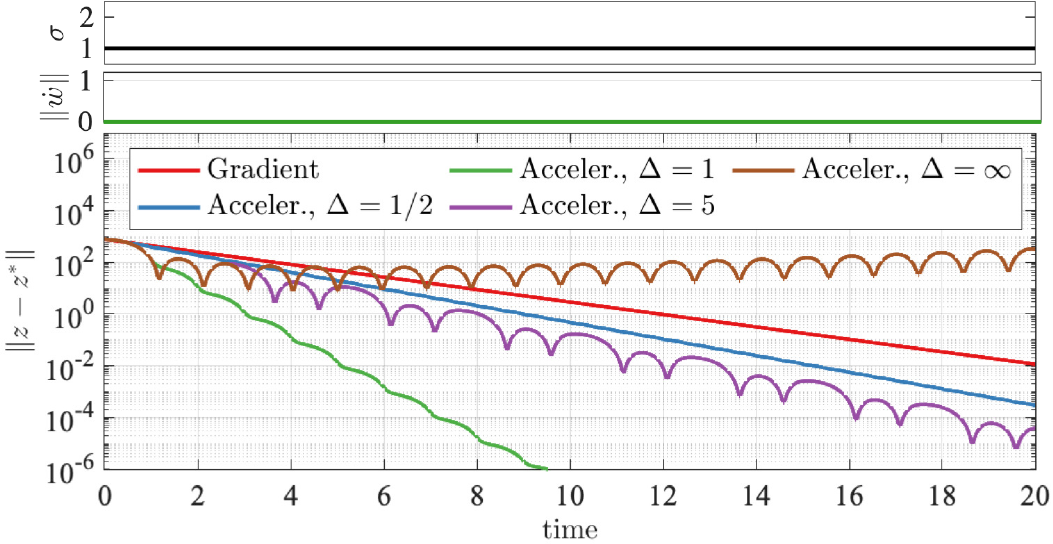}%
\caption{Comparison Gradient Flow vs Accelerated Gradient.}
\label{fig:thm4-2and3} 
\vspace{-0.15cm}
\end{figure}

\begin{figure}[tb]
\centering
\subfigure[]{\includegraphics[width=.48\columnwidth]{%
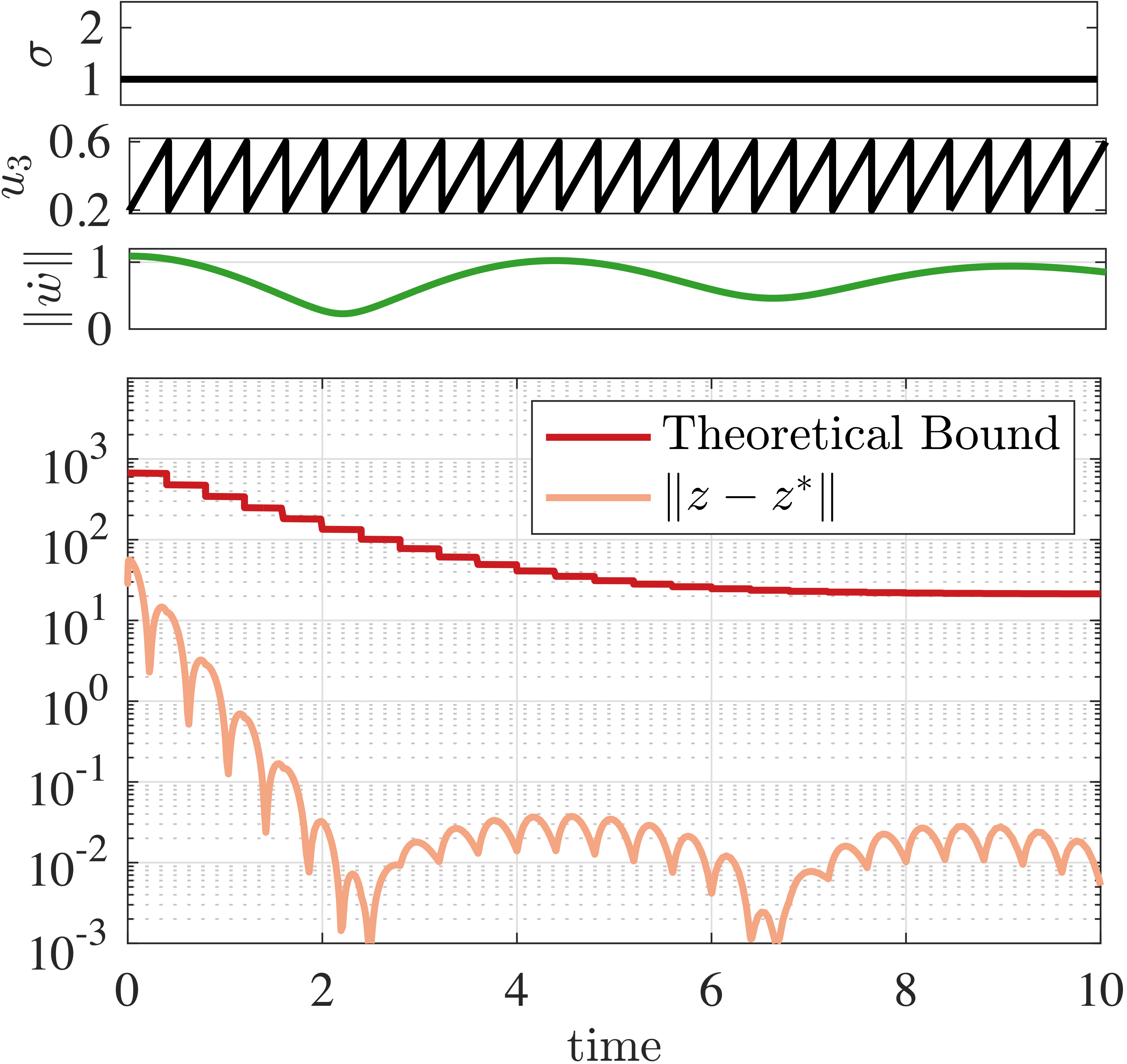}}%
\subfigure[]{\includegraphics[width=.48\columnwidth]{%
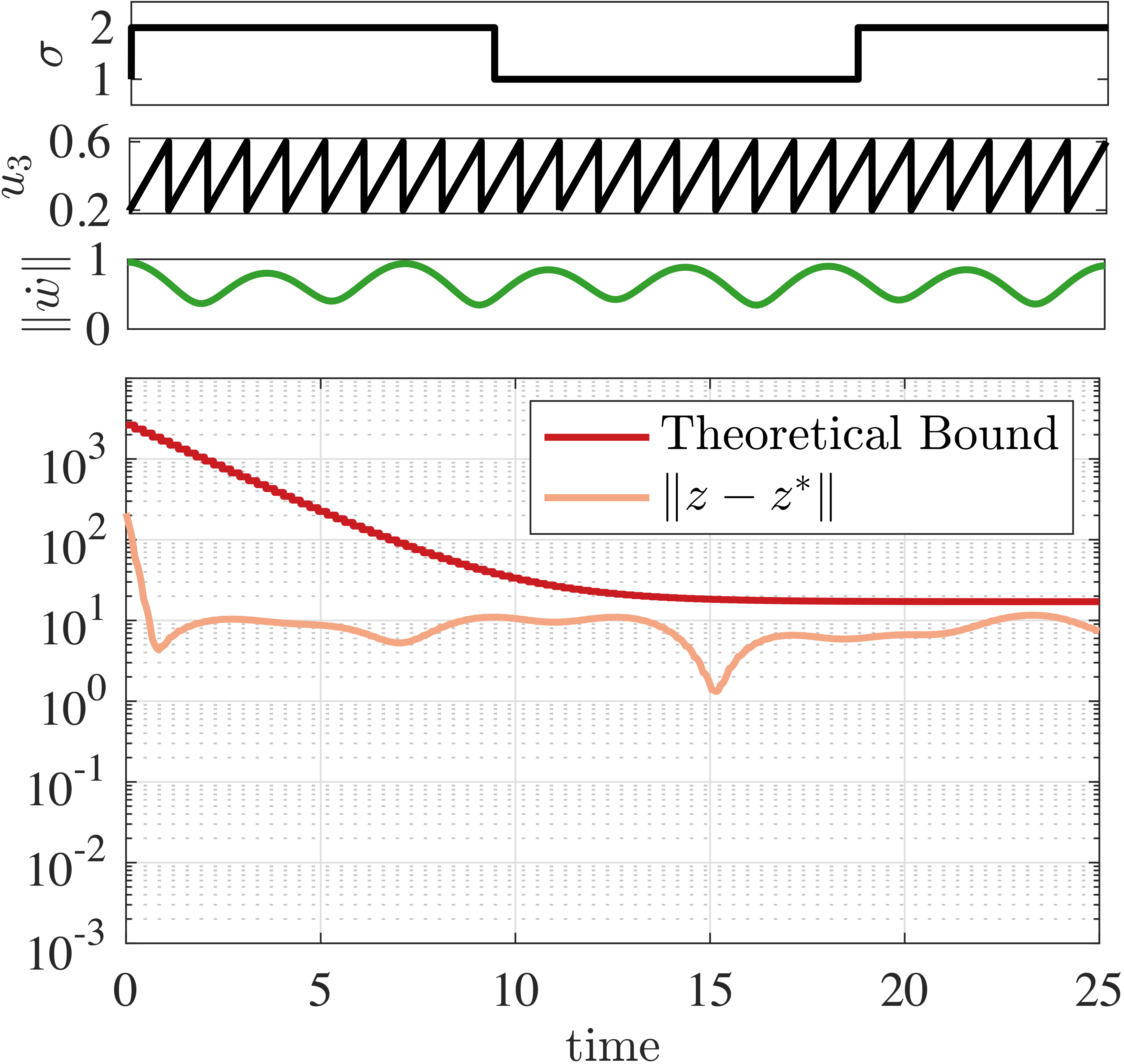}}
\caption{Tracking error of accelerated gradient controller. (a) No plant switching. 
(b) Switched plant.}
\label{fig:thm-4-4} 
\vspace{-0.1cm}
\end{figure}

\begin{figure}[tb]
\centering \subfigure[]{\includegraphics[width=.33\columnwidth]{%
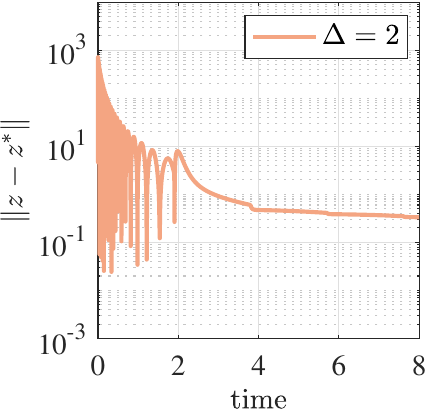}}%
\subfigure[]{\includegraphics[width=.33\columnwidth]{%
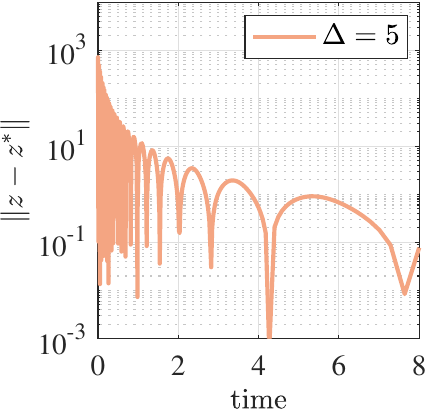}}%
\centering \subfigure[]{\includegraphics[width=.33\columnwidth]{%
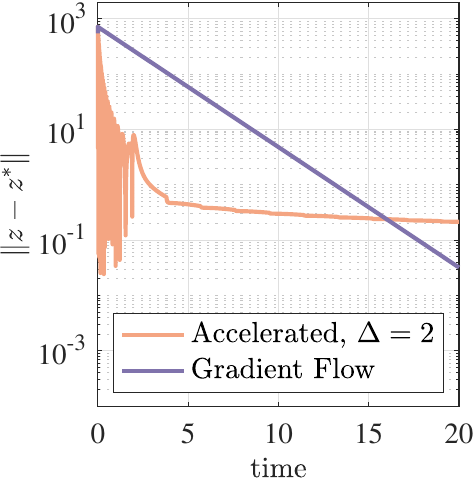}}
\caption{Tracking error of accelerated gradient controller with polynomial cost function. (a) $\Delta=2$. (b) $\Delta=5$. (c) Comparison with Gradient Flow.}
\label{fig:thm4-1} 
\end{figure}

To illustrate our results, we consider a plant with two modes 
$\mc S = \{1,2\}$, $n=10$ states, $m=5$ inputs, $p=5$ outputs, and $q=6$ 
exogenous disturbances. We consider cost functions
$\phi_u(u) = u^\tsp R u, \,\, \phi_y(y) = (y-\sbs{y}{ref})^\tsp Q (y-\sbs{y}{ref})$ 
where $R \in \real^{m \times m}$, $R \succ 0$, $Q^{p \times p}$, 
$Q \succ 0$, and $\sbs{y}{ref}\in \real^p$ is a constant reference 
signal. We first consider the gradient flow controller \eqref{eq:gradient_controller},  which generates the trajectories shown in Fig. \ref{fig:thm-3-3}. In particular,
Fig. \ref{fig:thm-3-3}-(a) shows the bound established in Theorem \ref{thm:EISS-Gradient}, when  the switching signal 
$\sigma$ is constant at all times. On the other hand, Fig. \ref{fig:thm-3-3}-(b) shows the bound corresponding to the case when the switching signal is time-varying. 
%
%
%
Next, we consider the hybrid accelerated gradient controller of Section \ref{sec:4}. When the disturbance $\omega_t$ is constant and the plant has a single operating mode (i.e. $\sigma=1$) Fig. \ref{fig:thm4-2and3} compares the performance of the hybrid controller versus the gradient-flow controller. \blue{Different reset parameters are considered, with $\Delta = 0.5, 1, 5$; we also consider $\Delta = \infty$, corresponding to \eqref{ODENesterov}.
In the latter case, it can be shown that the trajectories diverge.} Similarly, Fig. \ref{fig:thm-4-4}-(a) shows the bound of Theorem \ref{thm:EISS-Nesterov}, when $\omega_t$ is time-varying and  the switching signal $\sigma$ is constant at all times.  Fig. \ref{fig:thm-4-4}-(b) considers the case when the switching signal is  time-varying.
We note that in this simulation the time horizon has been rescaled in order to illustrate the behavior of the controller under switching of the plant.
%
Finally, Fig. \ref{fig:thm4-1} illustrates the bound in Theorem 
\ref{thm:EISS-Nesterov-Practical}. In this case, we consider a scalar plant ($n=m=p=1$) with a single mode  ($\vert \mc S \vert =1 $) and cost function $\phi_t(u)=  0.25(Gu+H\omega_t-\sbs{y}{ref})^4$, where $\sbs{y}{ref} \in \real$.
Note that $\phi_t$ satisfies Assumption \ref{ass:reverseLipschitz} and Assumption \ref{ass:lipschitzGradient} on compact sets.
Two important observations follow from Fig. \ref{fig:thm4-1}. 
First, the simulations suggest that in this case smaller restarting times $\Delta$ can improve transient performance.
Second, the simulation 
illustrates that the control signal converges only to a neighborhood of the optimal point, thus validating the conclusion of Theorem 
\ref{thm:EISS-Nesterov-Practical}.
The neighborhood of convergence can be characterized by observing that Assumption \ref{ass:reverseLipschitz} implies 
$(u-\sbs{u}{ref})^2 > \nu_0^2$, and thus $\nu_0^2 = \ell_0$.
This condition is illustrated at the top of Fig. \ref{fig:thm4-1}.
Finally, Fig. \ref{fig:thm4-1}-(c) compares the regulation error of the gradient flow controller versus the hybrid accelerated gradient controller. 
The figure shows that the accelerated gradient controller achieves faster convergence compared to the gradient descent-based controller, but the convergence is 
guaranteed only up to a neighborhood of the optimal points.

%



\vspace{-0.35cm}
\section{Conclusions}
\label{sec:7}
\vspace{-0.35cm}
We addressed the problem of online optimization of switched linear time invariant dynamical systems via  optimization-based controllers. We introduced two feedback controllers, one based on gradient descent flows, and a one based on a hybrid regularization of accelerated gradient systems with dynamic momentum. Under a suitable average dwell-time constraint on the switching signal, we established ISS properties for the closed-loop system with input being the time-derivative of the disturbance. This generalizes previous results on time-invariant optimization problems and non-switching plants. \blue{Future research directions will focus on developing tighter interconnection bounds between hybrid plants and hybrid controllers via small gain arguments.}


 \vspace{-0.35cm}
 \bibliographystyle{elsarticle-num}
 \bibliography{alias,bibliography}

\end{document}